\documentclass[a4paper]{amsart}%
\usepackage{amsfonts}
\usepackage{amssymb,latexsym}
\usepackage{amsmath}
\usepackage{amssymb}
\usepackage{graphicx}%
\theoremstyle{plain}
\newtheorem{theorem}{Theorem}
\newtheorem{lemma}[theorem]{Lemma}
\newtheorem{proposition}[theorem]{Proposition}
\theoremstyle{remark}
\newtheorem*{remark}{Remark}
\numberwithin{equation}{section}
\numberwithin{theorem}{section}
\newenvironment{caselist}
{
\begin{list}{\bf{Case \theenumi.}}
{
\usecounter{enumi}
\settowidth{\itemindent}{Case 3}
\addtolength{\itemindent}{\labelsep}
\settowidth{\labelwidth}{Case 3}
\setlength{\leftmargin}{0pt}
}
}
{\end{list}}

\begin{document}
\title[Almost associative operations generating a minimal clone]{Almost associative operations generating a minimal clone}
\author[T. Waldhauser]{Tam\'{a}s Waldhauser}
\address{Bolyai Institute\\
University of Szeged\\
Aradi v\'{e}rtan\'{u}k tere 1, H6720, Szeged, Hungary}
\email{twaldha@math.u-szeged.hu}
\thanks{{Research supported by the Hungarian National Foundation for Scientific
Research grant no. T 048809.}}
\date{}

\begin{abstract}
Characterizations of `almost associative' binary operations generating a
minimal clone are given for two interpretations of the term `almost
associative'. One of them uses the associative spectrum, the other one uses
the index of nonassociativity to measure how far an operation is from being associative.

\end{abstract}
\maketitle

\section{Minimal clones}

A \emph{clone} on a set $A$ is a set of finitary operations on $A$ that is
closed under composition of functions and contains all the projections. The
base set $A$ can be arbitrary; we will never assume finiteness in this paper.
If $\mathbb{A}=\left(  A;F\right)  $ is an algebra, then the set of term
functions, denoted by $Clo\left(  \mathbb{A}\right)  $, is a clone on $A$, the
\emph{clone of the algebra} $\mathbb{A}$. In this case $Clo\left(
\mathbb{A}\right)  $ is the smallest clone containing $F$, therefore we say
that $F$ \emph{generates} the clone, and we write $\left[  F\right]
=Clo\left(  \mathbb{A}\right)  $. (Clearly, every clone arises as the clone of
an algebra: we just need to pick a generating set for the clone, and let these
be the basic operations of the algebra.) We can also speak about $Clo\left(
\mathcal{V}\right)  $, the \emph{clone of a variety }$\mathcal{V}$. By this we
mean the clone of $\mathbb{F}_{\aleph_{0}}\left(  \mathcal{V}\right)  $, the
countably generated free algebra of $\mathcal{V}$.

The $n$-ary part of $Clo\left(  \mathbb{A}\right)  $, denoted by $Clo^{\left(
n\right)  }\left(  \mathbb{A}\right)  $ can be identified naturally with
$\mathbb{F}_{n}\left(  V\left(  \mathbb{A}\right)  \right)  $, the
$n$-generated free algebra of the variety generated by $\mathbb{A}$.
Projections correspond to variables under this identification: the first
binary projection is $e_{1}:\left(  x,y\right)  \mapsto x$ and the second one
is $e_{2}:\left(  x,y\right)  \mapsto y$, therefore we will sometimes think of
the variables $x$ and $y$ as projections (most of the time we will work with
binary operations).

All clones on a given set $A$ form a lattice with respect to inclusion; the
smallest element of this lattice is the \emph{trivial clone}, the clone of all
projections on $A$, while the greatest element is the clone of all finitary
operations on $A$. \emph{Minimal clones} are the atoms of this lattice, i.e. a
clone is minimal, if its only proper subclone is the trivial clone. A minimal
clone is generated by any of its nontrivial (i.e. non-projection) elements,
thus all minimal clones are singly generated, and therefore arise as clones of
algebras with just one basic operation. If $\mathbb{A=}\left(  A;f\right)  $
is such an algebra, then in order to prove that it has a minimal clone, one
needs to verify that $f\in\left[  g\right]  $ holds for every nontrivial $g\in
Clo\left(  \mathbb{A}\right)  $. This fact can be expressed by identities, so
if $\mathbb{A}$ has a minimal clone, then so does $V\left(  \mathbb{A}\right)
$, and if a variety $\mathcal{V}$ has a minimal clone, then the clone of any
algebra in $\mathcal{V}$ is either minimal or trivial.

To prove that a clone $\left[  f\right]  $ on $A$ is not minimal one needs to
find a nontrivial operation $g\in\left[  f\right]  $ such that $f\notin\left[
g\right]  $. This can be done for example by showing that there is an
equivalence relation $\rho$ on $A$ (a subset of $A$), such that $\rho$ is a
congruence (subuniverse) of the algebra $\left(  A;g\right)  $, but is not a
congruence (subuniverse) of $\left(  A;f\right)  $. (There is a general notion
of preservation of relations of arbitrary arity, and this gives a
Galois-correspondence between operations and relations on finite sets
\cite{BKKR Galois, G Galois}, but we will use only the previous observation,
which is valid for infinite sets as well.)

It is convenient to generate a minimal clone by a nontrivial operation of the
smallest arity. Minimal clones are classified with respect to this generator;
there are five types, and for two of them there is a complete characterization
of minimal clones (Rosenberg's Theorem, see \cite{R 5typ} and \cite{SzA clUA}).

One of the three types where the description of minimal clones is not complete
yet is the binary case. Clones of this type are generated by an idempotent
binary operation, so they can be (and will be) viewed as clones of idempotent
groupoids. (In this paper the term \emph{groupoid} refers to an algebra with a
single binary operation.) The basic operation of a groupoid will be denoted by
$f\left(  x,y\right)  =xy$, and by the \emph{dual} of $\mathbb{A=}\left(
A;f\right)  $ we mean the groupoid $\mathbb{A}^{d}=\left(  A;f^{d}\right)  $
with $f^{d}\left(  x,y\right)  =f\left(  y,x\right)  =yx$. Similarly,
$\mathcal{V}^{d}$ denotes the variety formed by the duals of the elements of
$\mathcal{V}$. Obviously, a groupoid has a minimal clone if and only if its
dual does (actually they have the very same clone).

A groupoid has a trivial clone if and only if it is a left or right zero
semigroup. The simplest examples of groupoids (or varieties) with a minimal
clone are semilattices and rectangular bands. Before giving more examples of
varieties with a minimal clone, let us make some notational conventions.

To save parentheses we use the notation $\overleftarrow{x_{1}\cdot\ldots\cdot
x_{n}}$ for the left-associated product $\left(  \cdots\left(  \left(
x_{1}x_{2}\right)  x_{3}\right)  \cdots\right)  x_{n}$, and similarly
$\overrightarrow{x_{1}\cdot\ldots\cdot x_{n}}$ for the right-associated
product $x_{1}\left(  \cdots\left(  x_{n-2}\left(  x_{n-1}x_{n}\right)
\right)  \cdots\right)  $. We abbreviate $\overleftarrow{x\cdot y\cdot
\ldots\cdot y} $ to $xy^{n}$ (where $n$ is certainly the number of $y$'s
appearing in the product). Analogously ${}^{n}\hspace{-1pt}xy$ stands for
$\overrightarrow{x\cdot\ldots\cdot x\cdot y}$.

Let $\mathcal{B}$ denote the variety defined by the identities $xx\approx
x,x\left(  xy\right)  \approx x\left(  yx\right)  \approx\left(  xy\right)
x\approx\left(  xy\right)  y\approx\left(  xy\right)  \left(  yx\right)
\approx xy$; let $\mathcal{C}_{p}$ be the variety of $p$\emph{-cyclic
groupoids} (cf. \cite{Plonka k-cyc}) defined by $xx\approx x,x\left(
yz\right)  \approx xy,\left(  xy\right)  z\approx\left(  xz\right)
y,xy^{p}\approx x$, and finally, let $\mathcal{D}$ be defined by $x\left(
yx\right)  \approx\left(  xy\right)  x\approx\left(  xy\right)  y\approx
\left(  xy\right)  \left(  yx\right)  \approx xy$ and $x\cdot\overleftarrow
{x\cdot y_{1}\cdot\ldots\cdot y_{n}}\approx x$ (for all $n\geq0$). The clone
of $\mathcal{B}$ and $\mathcal{D}$ is minimal, while the clone of
$\mathcal{C}_{p}$ is minimal iff $p$ is a prime. The minimality of the clone
of $\mathcal{B}$ and $\mathcal{D}$ is proved in \cite{LLPPP}; these are the
clones in parts (c) and (d) in Theorem 5.2. For the proof of the minimality of
the clones of $p$-cyclic groupoids see \cite{Plonka idred}. (From now on we
always assume that $p$ is a prime number, when we speak about $p$-cyclic groupoids.)

The following propositions show the usefulness of \emph{absorption identities}
in the study of minimal clones. These are identities of the form $t\approx x$,
i.e. identities with a single variable on one side. The proofs of these
propositions can be found in \cite{LLPPP} and \cite{KK}.

\begin{proposition}
\label{prop absorption id.} Let $\mathcal{V}$ be a variety with a minimal
clone, and let $\mathbb{A\in}\mathcal{V}$ have a nontrivial (hence minimal)
clone. Then $\mathcal{V}$ satisfies every absorption identity that holds in
$\mathbb{A}$.
\end{proposition}

\begin{proof}
See Lemma 2.1 in \cite{LLPPP} or Lemma 3.6 in \cite{KK}.
\end{proof}

\begin{proposition}
\label{prop p-cyclic, rect. band}Let $\mathcal{V} $ be a variety with a
minimal clone, and suppose that $\mathcal{V}$ contains a $p$-cyclic groupoid
(rectangular band) with a nontrivial clone. Then $\mathcal{V}$ is the variety
of $p$-cyclic groupoids (rectangular bands).
\end{proposition}

\begin{proof}
By the previous proposition, it suffices to show that $p$-cyclic groupoids and
rectangular bands are axiomatizable by absorption identities. For $p$-cyclic
groupoids such an axiomatization is given in Lemma 3.10 of \cite{KK}, and the
method described in Lemma 2.3 of \cite{LLPPP} yields (almost) the same
identities. For rectangular bands see Lemma 3.8 of \cite{KK} or Theorem 5.2
(b) of \cite{LLPPP} for a list of absorption identities. Note that the
conclusion of the proposition says that \textquotedblleft$\mathcal{V}$ is
\emph{the} variety of $p$-cyclic groupoids (rectangular
bands)\textquotedblright, not that \textquotedblleft$\mathcal{V}$ is \emph{a}
variety of $p$-cyclic groupoids (rectangular bands)\textquotedblright. This is
because the only nontrivial subvariety of $\mathcal{C}_{p}$ is the variety of
left zero semigroups (see the last paragraph of the proof of Lemma 3.5 in
\cite{KK} or Corollary 2.1 in \cite{LLPPP}), and clearly the variety of
rectangular bands does not have a proper subvariety with a nontrivial clone either.
\end{proof}

\begin{proposition}
\label{prop D}Let $\mathcal{V}$ be a variety with a minimal clone satisfying
the identities $xx\approx x,x\left(  yx\right)  \approx\left(  xy\right)
x\approx\left(  xy\right)  y\approx\left(  xy\right)  \left(  yx\right)
\approx xy,x\left(  xy\right)  \approx x$. Then $\mathcal{V}$ is a subvariety
of $\mathcal{D}$.
\end{proposition}

\begin{proof}
This is part (d) of Theorem 5.2 in \cite{LLPPP}. The identities listed here
are sufficient to determine the two-generated free algebra of $\mathcal{V}$.
Its multiplication table is the following (the four elements have to be
distinct, since otherwise $Clo\left(  \mathcal{V}\right)  $ would be
trivial).
\[%
\begin{tabular}
[c]{l|llll}%
$\cdot$ & $x$ & $y$ & $xy$ & $yx$\\\hline
$x$ & $x$ & $xy$ & $x$ & $xy$\\
$y$ & $yx$ & $y$ & $yx$ & $y$\\
$xy$ & $xy$ & $xy$ & $xy$ & $xy$\\
$yx$ & $yx$ & $yx$ & $yx$ & $yx$%
\end{tabular}
\
\]

It is not hard to check that this groupoid satisfies every identity of the
form $x\cdot\overleftarrow{x\cdot y_{1}\cdot\ldots\cdot y_{n}}\approx x$ (this
is a special case of Lemma 4.2 in \cite{LLPPP}). These are absorption
identities, therefore we can apply Proposition \ref{prop absorption id.} with
$\mathbb{A=F}_{2}\left(  \mathcal{V}\right)  $ to show that $\mathcal{V}$
satisfies these identities, too. The remaining identities in the definition of
$\mathcal{D}$ are the same as that were assumed.
\end{proof}

Characterizing minimal clones in general is a hard task, even in the binary
case. All known results describe minimal clones under certain restrictions
\cite{Cs 3all, Cs cons, JQ cons, KK, KSz comm, LLPPP, Post, Szcz, W 4maj, W
wa}. Another result of this kind is the description of associative binary
operations generating a minimal clone \cite{PPP prep, SzM minfcs}: a semigroup
has a minimal clone iff it is a rectangular band, a left regular band
(idempotent semigroup satisfying $xyx\approx xy$) or a right regular band
(dual of a left regular band). Note that left and right regular bands belong
to the varieties $\mathcal{B}$ and $\mathcal{B}^{d}$, respectively. In this
paper we slightly generalize this result by characterizing `almost
associative' binary operations generating a minimal clone. To explain what we
mean by being `almost associative', we need a way to measure how far a certain
operation is from being associative. We discuss two such measures: the
associative spectrum and the index of nonassociativity. In Section
\ref{Section spec} we characterize groupoids with a minimal clone and small
associative spectrum (Theorem \ref{THM spec}), and in Section \ref{Section SH}
we describe groupoids with a minimal clone and small index of nonassociativity
(Theorem \ref{THM SH}).

\section{Minimal clones with small associative spectrum\label{Section spec}}

One way of measuring associativity is possible by considering the identities
implied by associativity, and somehow counting how many of these are (not)
satisfied. To make this more precise, let us say that $B$ is a
\emph{bracketing}, if $B$ is a groupoid term, and each variable occurs exactly
once in $B$. If these variables are $x_{1},x_{2},\ldots,x_{n}$ and they appear
in this order (as we will suppose most of the time), then $B$ is nothing else
but a way to put brackets into the product $x_{1}\cdot\ldots\cdot x_{n}$ such
that the order of the $n-1$ multiplications is well determined. In this case
we say that $B$ is a bracketing of the product $x_{1}\cdot\ldots\cdot x_{n}$,
and we write $B=B\left(  x_{1},\ldots,x_{n}\right)  $. The number of variables
appearing in $B$ is called the size of $B$, and is denoted by $\left\vert
B\right\vert $.

In every bracketing there is an outermost multiplication, and this splits the
bracketing into two parts, the \emph{left factor} and the \emph{right factor}
of the bracketing. Let $B=B\left(  x_{1},\ldots,x_{n}\right)  $, and let $P,Q$
be the left and right factors of $B$. Then $B=PQ$, and $P=P\left(
x_{1},\ldots,x_{k}\right)  ,Q=Q\left(  x_{k+1},\ldots,x_{n}\right)  $, where
$k=\left\vert P\right\vert $. Sometimes we will use the notation $l\left(
B\right)  $ for the left factor of $B$.

The number of bracketings of the product $x_{1}\cdot\ldots\cdot x_{n}$ is
$C_{n-1}=\frac{1}{n}\binom{2n-2}{n-1}$, the $\left(  n-1\right)  $st Catalan
number. In a semigroup, all of these $C_{n-1}$ many terms induce the same term
function, but in an arbitrary groupoid they may induce more than one term
function. Intuitively, the more term functions of this kind there are, the
less associative the multiplication is. Therefore we define the
\emph{associative spectrum} of a groupoid $\mathbb{A}$ to be the sequence
$s_{\mathbb{A}}\left(  1\right)  ,s_{\mathbb{A}}\left(  2\right)
,\ldots,s_{\mathbb{A}}\left(  n\right)  ,\ldots$, where $s_{\mathbb{A}}\left(
n\right)  $ is the number of different term functions on $\mathbb{A}$ arising
from bracketings of $x_{1}\cdot\ldots\cdot x_{n}$. Thus the associative
spectrum gives (only quantitative) information about identities of the form
$B_{1}\left(  x_{1},\ldots,x_{n}\right)  \approx B_{2}\left(  x_{1}%
,\ldots,x_{n}\right)  $ satisfied by the groupoid. The associative spectrum
was introduced and investigated in \cite{CsW spec}.

Clearly, $s_{\mathbb{A}}\left(  1\right)  =s_{\mathbb{A}}\left(  2\right)  =1$
for every groupoid $\mathbb{A}$, and $s_{\mathbb{A}}\left(  3\right)  =1$ iff
$\mathbb{A}$ is a semigroup. In the latter case $s_{\mathbb{A}}\left(
n\right)  =1$ for all $n$ by the general law of associativity. The smallest
possible spectrum for a nonassociative multiplication is $1,1,2,1,1,\ldots$,
so we could say that a binary operation is almost associative, if its spectrum
is this sequence. However, there is no groupoid having a minimal clone with
this spectrum (not even an idempotent groupoid) as we will see later.
Therefore we have to be more generous: in Theorem \ref{THM spec} we determine
groupoids with a minimal clone satisfying $s\left(  4\right)  <5=C_{3}$. First
we prove three theorems which show that certain identities of the form
$B_{1}\left(  x_{1},\ldots,x_{n}\right)  \approx B_{2}\left(  x_{1}%
,\ldots,x_{n}\right)  $ cannot hold in nonassociative groupoids with a minimal
clone, and then we discuss the four-variable case in detail. (In the first two
theorems we actually assume only idempotence.)

\begin{theorem}
\label{thm 1-4}If an idempotent groupoid satisfies the identity
\begin{equation}
\overleftarrow{x_{1}\cdot\ldots\cdot x_{n}}\approx\overrightarrow{x_{1}%
\cdot\ldots\cdot x_{n}} \label{bal=jobb}%
\end{equation}
for some $n\geq3$, then it is a semigroup.
\end{theorem}

\begin{proof}
Applying (\ref{bal=jobb}) with $x_{1}=\ldots=x_{k}=x,x_{k+1}=\ldots=x_{n}=y $
we obtain%
\begin{equation}
xy^{n-k}\approx\overleftarrow{x\cdot\ldots\cdot x\cdot y\cdot\ldots\cdot
y}\approx\overrightarrow{x\cdot\ldots\cdot x\cdot y\cdot\ldots\cdot y}%
\approx{}^{k}\hspace{-1pt}xy \label{hatvanyok}%
\end{equation}
for $1\leq k\leq n-1.$

Let us use (\ref{bal=jobb}) again, for $x_{1}=x,x_{2}=u=xy^{2}\approx{}%
^{n-2}\hspace{-1pt}xy,x_{3}=\ldots=x_{n}=y$:%
\begin{equation}
\left(  xu\right)  y^{n-2}\approx\overleftarrow{x\cdot u\cdot y\cdot
\ldots\cdot y}\approx\overrightarrow{x\cdot u\cdot y\cdot\ldots\cdot y}\approx
x\left(  uy\right)  . \label{ize}%
\end{equation}
The left hand side is $\left(  xu\right)  y^{n-2}\approx{}\left(
^{n-1}\hspace{-1pt}xy\right)  y^{n-2}\approx\left(  xy\right)  y^{n-2}\approx
xy^{n-1}\approx xy$ (we used (\ref{hatvanyok}) twice, with $k=n-1$ and $k=1$
respectively). We can compute the right hand side of (\ref{ize}) in a similar
manner: $x\left(  uy\right)  \approx x\left(  xy^{3}\right)  \approx x\left(
{}^{n-3}\hspace{-1pt}xy\right)  \approx{}^{n-2}\hspace{-1pt}xy\approx xy^{2}$.
Thus we have $xy\approx xy^{2}$, i.e. right multiplications are idempotent.

Finally, to prove associativity, we write up (\ref{bal=jobb}) one more time:%
\[
\left(  xy\right)  z^{n-2}\approx\overleftarrow{x\cdot y\cdot z\cdot
\ldots\cdot z}\approx\overrightarrow{x\cdot y\cdot z\cdot\ldots\cdot z}\approx
x\left(  yz\right)  .
\]
By the idempotence of right multiplication (by $z$) the left hand side reduces
to $\left(  xy\right)  z$, and therefore associativity is established.
\end{proof}

\begin{theorem}
\label{thm 1-2 ,4-5}An idempotent groupoid satisfying the following two
identities for some $n\geq3$, must be a semigroup.%
\begin{align*}
x_{0}\cdot\overleftarrow{x_{1}\cdot\ldots\cdot x_{n}}  &  \approx x_{0}%
\cdot\overrightarrow{x_{1}\cdot\ldots\cdot x_{n}}\\
\overleftarrow{x_{1}\cdot\ldots\cdot x_{n}}\cdot x_{0}  &  \approx
\overrightarrow{x_{1}\cdot\ldots\cdot x_{n}}\cdot x_{0}%
\end{align*}

\end{theorem}

\begin{proof}
Substituting $\overleftarrow{x_{1}\cdot\ldots\cdot x_{n}}$ into $x_{0}$ in the
first identity we have%
\[
\overleftarrow{x_{1}\cdot\ldots\cdot x_{n}}\approx\overleftarrow{x_{1}%
\cdot\ldots\cdot x_{n}}\cdot\overrightarrow{x_{1}\cdot\ldots\cdot x_{n}}.
\]
by idempotence. Similarly, if we substitute $\overrightarrow{x_{1}\cdot
\ldots\cdot x_{n}}$ for $x_{0}$ in the second identity, then we get
\[
\overleftarrow{x_{1}\cdot\ldots\cdot x_{n}}\cdot\overrightarrow{x_{1}%
\cdot\ldots\cdot x_{n}}\approx\overrightarrow{x_{1}\cdot\ldots\cdot x_{n}},
\]
and thus (\ref{bal=jobb}), hence also associativity follows by the previous theorem.
\end{proof}

\begin{theorem}
\label{thm 1-5}If a groupoid has a minimal clone and satisfies%
\begin{equation}
\overleftarrow{x_{1}\cdot\ldots\cdot x_{n}}\approx x_{1}\cdot\overleftarrow
{x_{2}\cdot\ldots\cdot x_{n}} \label{nulla}%
\end{equation}
for some $n\geq3$, then it is a semigroup.
\end{theorem}

\begin{proof}
The case $n=3$ is trivial, so let us suppose that $n\geq4$. First we draw a
consequence of (\ref{nulla}) and idempotence (putting $x$ and $z$ for $x_{1}$
and $x_{n}$, and $y$ for the rest of the variables):%
\begin{equation}
\left(  xy^{n-2}\right)  z\approx x\left(  yz\right)  . \label{ketto}%
\end{equation}
As a special case (with $z=y$) we get%
\begin{equation}
xy^{n-1}\approx xy. \label{egy}%
\end{equation}

Now we suppose that $\mathbb{A=}\left(  A;\cdot\right)  $ is a groupoid with a
minimal clone that satisfies identity (\ref{nulla}). The binary operation
$s\left(  x,y\right)  =xy^{n-2}$ belongs to the clone of $\mathbb{A}$,
therefore if it is nontrivial, then $\left[  s\right]  $ contains the basic
operation $f\left(  x,y\right)  =xy$.

Suppose that $a$ and $b$ are arbitrary elements of $A$ such that
$c=(ab)a^{n-3}\neq a.$ We claim that $s$ is a semilattice operation on the
two-element set $\left\{  a,c\right\}  $. With the help of (\ref{egy}) we see
that $s\left(  c,a\right)  =\left(  (ab)a^{n-3}\right)  a^{n-2}=(ab)a^{2n-5}%
=\left(  (ab)a^{n-1}\right)  a^{n-4}=\left(  \left(  ab\right)  a\right)
a^{n-4}=(ab)a^{n-3}=c$. To compute $s\left(  a,c\right)  $ let us first
consider $ac$:%
\begin{equation}
ac=a\left(  (ab)a^{n-3}\right)  =\left(  (aa)b\right)  a^{n-3}=(ab)a^{n-3}=c.
\label{a(aba^n-3)}%
\end{equation}
In the middle two steps we used identity (\ref{nulla}) and idempotence. Now it
is easy to conclude that $s\left(  a,c\right)  =ac^{n-2}=c$, proving that $s $
is indeed a semilattice operation on $\left\{  a,c\right\}  $.

Since $f\in\left[  s\right]  $, the restriction of $f$ to $\left\{
a,c\right\}  $ is either trivial, or coincides with $s.$ In the latter case we
have $f\left(  c,a\right)  =c$, so%
\begin{equation}
\left(  (ab)a^{n-3}\right)  a=(ab)a^{n-2}=(ab)a^{n-3}. \label{pre-harom}%
\end{equation}

If $f$ is trivial on our two-element set, then it has to be a second
projection, because $f\left(  a,c\right)  =ac=c$ as we have already observed
in (\ref{a(aba^n-3)}). Thus we have $f\left(  c,a\right)  =ca=a$, which means
that $(ab)a^{n-2}=a$. Multiplying by $a$ from the right we get $(ab)a^{n-1}=a
$, therefore $\left(  ab\right)  a=a$ by (\ref{egy}). If we multiply both
sides of this equality $n-4$ times by $a$, then we get $(ab)a^{n-3}=a$, i.e.
$c=a$, contrary to our assumption.

If $(ab)a^{n-3}=a$ holds for $a,b\in A$, then (\ref{pre-harom}) holds
trivially. Thus we have proved that if a groupoid $\mathbb{A}$ has a minimal
clone, and satisfies (\ref{nulla}), then (\ref{pre-harom}) holds for all
$a,b\in A$. In other words, $\mathbb{A}$ satisfies the following identity.%
\begin{equation}
(xy)x^{n-3}\approx(xy)x^{n-2} \label{harom}%
\end{equation}

It suffices to show now that (\ref{nulla}) and (\ref{harom}) together with
idempotence imply associativity. Let us multiply both sides of (\ref{harom})
by $x$ from the right. We get $(xy)x^{n-2}\approx(xy)x^{n-1}$ and then
(\ref{egy}) shows that $(xy)x^{n-2}\approx(xy)x$. Therefore $\left(
(xy)x^{n-2}\right)  z\approx\left(  (xy)x\right)  z$ also holds. The left hand
side of this identity reduces to $\left(  xy\right)  \left(  xz\right)  $
according to (\ref{ketto}), with $xy$, $x$ and $z$ playing the role of $x$,
$y$ and $z$, respectively. Thus we have obtained the following identity.%
\begin{equation}
\left(  \left(  xy\right)  x\right)  z\approx\left(  xy\right)  \left(
xz\right)  \label{negy}%
\end{equation}

Now we go back to (\ref{harom}), and this time we multiply it by $y$ from the
left. The left hand side becomes $y\left(  (xy)x^{n-3}\right)  $, which turns
to $\left(  \left(  yx\right)  y\right)  x^{n-3}$ if we apply (\ref{nulla}).
With the help of (\ref{negy}) and idempotence we can simplify this expression:
$\left(  \left(  yx\right)  y\right)  x^{n-3}\approx\left(  \left(  \left(
yx\right)  y\right)  x\right)  x^{n-4}\approx\left(  \left(  yx\right)
\left(  yx\right)  \right)  x^{n-4}\approx\left(  yx\right)  x^{n-4}\approx
yx^{n-3}$. The right hand side of (\ref{harom}) becomes $y\left(  \left(
xy\right)  x^{n-2}\right)  $. This can be considered as a product of $n$
factors, if we keep the $x$ and the $y$ in the middle together. We can
rearrange this product according to (\ref{nulla}), and we get $\left(
y\left(  xy\right)  \right)  x^{n-2}$. The $y\left(  xy\right)  $ at the
beginning of this term can be written as $y\cdot\overleftarrow{x\cdot
\ldots\cdot x\cdot y}$, and an application of (\ref{nulla}) yields
$\overleftarrow{y\cdot x\cdot\ldots\cdot x\cdot y}\approx\left(
yx^{n-2}\right)  y$. Substituting this back into the original expression we
get $\left(  y\left(  xy\right)  \right)  x^{n-2}\approx\left(  \left(
yx^{n-2}\right)  y\right)  x^{n-2}$. If we consider $yx^{n-2}$ as one factor,
then this is again a (left-associated) product of $n$ factors, and we can use
(\ref{nulla}) one more time: $\left(  \left(  yx^{n-2}\right)  y\right)
x^{n-2}\approx\left(  yx^{n-2}\right)  \left(  yx^{n-2}\right)  $. Clearly
this is just $yx^{n-2}$, and if we compare the results we have obtained from
the two sides of (\ref{harom}) we can conclude the following identity.%
\[
yx^{n-3}\approx yx^{n-2}%
\]

Multiplying this by $x$ we get $yx^{n-2}\approx yx^{n-1}\approx yx$ by
(\ref{egy}). Now the left hand side of (\ref{ketto}) can be simplified as
$\left(  xy^{n-2}\right)  z\approx\left(  xy\right)  z$, and therefore
associativity follows.
\end{proof}

\begin{remark}
Idempotence and identity (\ref{nulla}) for $n\geq4$ do not imply
associativity, as we can see from the following example. For every $k\geq2$ we
define a groupoid $\mathbb{A}_{k}$ on the set $A_{k}=\mathbb{Z}\hspace
{0cm}_{k}\dot{\cup}\left\{  e\right\}  $ by%
\[
xy=%
\begin{cases}
y & \text{if }y\neq e;\\
x+1 & \text{if }y=e\neq x;\\
e & \text{if }y=e=x.
\end{cases}
\]

This groupoid is idempotent, but not associative, because $\left(  0e\right)
e=2\neq1=0\left(  ee\right)  $. Let $B\left(  x_{1},\ldots,x_{n}\right)  $ be
a bracketing, and let $l_{i}$ denote the left depth of $x_{i}$ in $B$ (see
\cite{CsW spec} for the definition of left depth). It is not hard to prove by
induction on $n$, that for any $c_{1},\ldots,c_{n}\in A_{k}$ we have $B\left(
c_{1},\ldots,c_{n}\right)  =c_{i}+l_{i}$ if $c_{i}$ is the last element of the
sequence $c_{1},\ldots,c_{n}$ that is different from $e$ (if there is no such
element, then clearly $B\left(  c_{1},\ldots,c_{n}\right)  =e $). Thus two
bracketings give the same term function on $\mathbb{A}_{k}$ iff their left
depth sequences are congruent modulo $k.$ The left depth sequence of the
bracketing on the left hand side of (\ref{nulla}) is $\left(
n-1,n-2,n-3,\ldots,1,0\right)  $ and that of the right hand side is $\left(
1,n-2,n-3,\ldots,1,0\right)  $. Hence $\mathbb{A}_{k}$ satisfies (\ref{nulla})
iff $k$ divides $n-2$. For example, $\mathbb{A}_{n-2}$ is an idempotent
nonassociative groupoid satisfying (\ref{nulla}).

The associative spectrum of $\mathbb{A}_{k}$ is the same as that of the
operation $x+\varepsilon y$ on $\mathbb{C}\hspace{0cm}$, where $\varepsilon$
is a primitive $k$-th root of unity: both count the number of \emph{zag
sequences} modulo $k$ (cf. \cite{CsW spec} 2.8. and 6.4.). If $k=2$, then we
have $\varepsilon=-1$, and the spectrum is $2^{n-2}$ (cf. \cite{CsW spec}
3.1.). For $k=3$ the spectrum is sequence A005773 in the Encyclopedia
\cite{enc}; this sequence is related to Motzkin numbers (A001006). The
spectrum for $k=4$ does not appear in the Encyclopedia, but the superseeker
found that it is a transform of the sequence A036765
\end{remark}

Let us now turn to the investigation of four-variable `associativity
conditions'. There are five bracketings of size four:%
\begin{align*}
B_{1}  &  =x\left(  y\left(  zu\right)  \right)  ;\\
B_{2}  &  =x\left(  \left(  yz\right)  u\right)  ;\\
B_{3}  &  =\left(  xy\right)  \left(  zu\right)  ;\\
B_{4}  &  =\left(  \left(  xy\right)  z\right)  u;\\
B_{5}  &  =\left(  x\left(  yz\right)  \right)  u.
\end{align*}
Many of the possible $\binom{5}{2}$ identities cannot be satisfied by a
nonassociative idempotent groupoid. For example, identifying $z$ and $u$ in
$B_{1}$ and $B_{3}$ we see that $B_{1}\approx B_{3}$ implies associativity if
idempotence is assumed. A similar argument works for $B_{3}\approx B_{4}$ and
$B_{2}\approx B_{5}$. For $B_{2}\approx B_{3}$ we need two steps: multiplying
both sides by a variable from the left yields $x\left(  y\left(  \left(
zu\right)  v\right)  \right)  \approx x\left(  \left(  yz\right)  \left(
uv\right)  \right)  $ (after renaming the variables), while replacing $u$ with
$uv$ gives $x\left(  \left(  yz\right)  \left(  uv\right)  \right)
\approx\left(  xy\right)  \left(  z\left(  uv\right)  \right)  $. Now
$x\left(  y\left(  \left(  zu\right)  v\right)  \right)  \approx\left(
xy\right)  \left(  z\left(  uv\right)  \right)  $ follows by transitivity, and
identifying $z,u$ and $v$ we get $x\left(  yz\right)  \approx\left(
xy\right)  z$. We can treat $B_{3}\approx B_{5}$ similarly (this is actually
the dual of $B_{2}\approx B_{3}$).

Specializing Theorems \ref{thm 1-4} and \ref{thm 1-5} to $n=4$ we see that
$B_{1}\approx B_{4}$ and $B_{2}\approx B_{4}$ cannot hold in a nonassociative
groupoid with a minimal clone, and neither can $B_{1}\approx B_{5}$, because
it is the dual of $B_{2}\approx B_{4}$. Only three possibilities remain: our
groupoid satisfies $B_{1}\approx B_{2}$ or $B_{4}\approx B_{5}$ or both.
Theorem \ref{thm 1-2 ,4-5} shows that the third case is impossible, hence we
can conclude that if a groupoid $\mathbb{A}$ has a minimal clone, and
$1<s_{\mathbb{A}}\left(  4\right)  <5$ holds for its spectrum, then
$s_{\mathbb{A}}\left(  4\right)  =4$, and $\mathbb{A}$ satisfies either
$B_{1}\approx B_{2}$ or its dual, but not both. We are going to characterize
such groupoids in the next theorem, but first we need three lemmas. Let
$\mathcal{A}$ denote the variety defined by $B_{1}\approx B_{2}$, i.e.
$x\left(  y\left(  zu\right)  \right)  \approx x\left(  \left(  yz\right)
u\right)  $.

\begin{lemma}
\label{lemma general ass}If $t_{1}\approx t_{2}$ is an identity that is true
in every semigroup, then $\mathcal{A}$ satisfies $xt_{1}\approx xt_{2}$ (where
$x$ is an arbitrary variable).
\end{lemma}

\begin{proof}
If $t_{1}\approx t_{2}$ holds in the variety of semigroups, then $t_{1}$ and
$t_{2}$ are two bracketings of the same product. Therefore it suffices to
prove that $\mathcal{A}$ satisfies $x\cdot B\left(  x_{1},\ldots,x_{n}\right)
\approx x\cdot\overrightarrow{x_{1}\cdot\ldots\cdot x_{n}}$ for any bracketing
$B\left(  x_{1},\ldots,x_{n}\right)  $. We prove this by induction on $n$.

Repeatedly applying $x\left(  \left(  yz\right)  u\right)  \approx x\left(
y\left(  zu\right)  \right)  $ we can transform $x\cdot B\left(  x_{1}%
,\ldots,x_{n}\right)  $ to the form $x\cdot\left(  x_{1}\cdot B^{\prime
}\left(  x_{2},\ldots,x_{n}\right)  \right)  $. By the induction hypothesis we
have that $x_{1}\cdot B^{\prime}\left(  x_{2},\ldots,x_{n}\right)  \approx
x_{1}\cdot\overrightarrow{x_{2}\cdot\ldots\cdot x_{n}}=\overrightarrow
{x_{1}\cdot\ldots\cdot x_{n}}$ holds in $\mathcal{A}$, hence $x\cdot B\left(
x_{1},\ldots,x_{n}\right)  \approx x\cdot\overrightarrow{x_{1}\cdot\ldots\cdot
x_{n}}$ is true as well. (Note that we did nothing else but gave a proof for
the general law of associativity, but we had to avoid implications of the form
$p\approx q\Rightarrow pr\approx qr$).
\end{proof}

\begin{lemma}
\label{lemma v-w}Let $\mathcal{V}$ be a subvariety of $\mathcal{A}$, and let
$\mathcal{W}$ be the intersection of $\mathcal{V}$ and the variety of
semigroups. If an identity $t_{1}\approx t_{2}$ holds in $\mathcal{W}$, then
$xt_{1}\approx xt_{2}$ holds in $\mathcal{V}$ (where $x$ is an arbitrary variable).
\end{lemma}

\begin{proof}
Let $\Theta_{\mathcal{V}},\Theta_{\mathcal{W}},\Theta_{sgr}$ denote the
equational theories of $\mathcal{V},\mathcal{W}$ and the variety of
semigroups, respectively. These are fully invariant congruences of the free
groupoid on countably many generators, and $\Theta_{\mathcal{W}}$ equals
$\Theta_{\mathcal{V}}\vee\Theta_{sgr}$, i.e. the transitive closure of
$\Theta_{\mathcal{V}}\cup\Theta_{sgr}$. Therefore, if $\mathcal{W}$ satisfies
an identity $t_{1}\approx t_{2}$, then there are terms $p_{1},\ldots,p_{n}$
such that $p_{1}=t_{1}$, $p_{n}=t_{2}$ and $p_{i}\approx p_{i+1}$ holds in
$\mathcal{V}$ if $i$ is odd, and $p_{i}\approx p_{i+1}$ is a semigroup
identity if $i$ is even. Then $xp_{i}\approx xp_{i+1}$ is true in
$\mathcal{V}$ for every $i$ and any variable $x$. (For odd $i$'s this is
obvious; for even ones it follows from the previous lemma.) Now $xt_{1}\approx
xt_{2}$ follows by transitivity.
\end{proof}

The next lemma is based on the method used in the proof of Lemma 3.8 in
\cite{KSz comm}, and is basically just a slight generalization of the
situation considered there.

\begin{lemma}
\label{lemma minimono}Suppose that $\mathbb{A}$ is a groupoid with a minimal
clone, and $M$ is a subset of $Clo^{\left(  2\right)  }\left(  \mathbb{A}%
\right)  $ containing the first projection and at least one nontrivial
element, such that for all $f,g,h\in M$%
%TCIMACRO{\TeXButton{roman}{\renewcommand{\theenumi}{(\roman{enumi})}
%\renewcommand{\labelenumi}{\theenumi}}}%
%BeginExpansion
\renewcommand{\theenumi}{(\roman{enumi})}
\renewcommand{\labelenumi}{\theenumi}%
%EndExpansion

\begin{enumerate}
\item \label{lemma minimono i}$f\left(  g,h\right)  =g$

\item \label{lemma minimono ii}$f\left(  g,h^{d}\right)  =f\left(
g,e_{2}\right)  \in M.$
\end{enumerate}

\noindent Then $\mathbb{A}$ or its dual belongs to the variety $\mathcal{D}$
or $\mathcal{C}_{p}$ for some prime number $p$.
\end{lemma}

\begin{proof}
Let us recall that $e_{1}$ and $e_{2}$ are the first and second binary
projection respectively (we can write $g^{d}$ as $g\left(  e_{2},e_{1}\right)
$ with this notation). Note that $e_{2}=e_{1}^{d}$, hence
\ref{lemma minimono ii} means that $f\left(  g,h^{d}\right)  $ does not depend
on $h$ (as long as $h\in M$). We have $e_{1}\in M$, but $e_{2}\in$ $M$ is
impossible, because then \ref{lemma minimono ii} would imply (with $f=e_{2}$)
that $h^{d}=e_{2}$ for every $h\in M$, contradicting that $M$ has at least two
elements. If $f\in M$ is nontrivial and $f^{d}$ also belongs to $M$, then we
have $f\left(  e_{1},f^{d}\right)  =e_{1}$ by \ref{lemma minimono i}, and
$f\left(  e_{1},f^{d}\right)  =f\left(  e_{1},e_{2}\right)  =f$ by
\ref{lemma minimono ii}, hence $f=e_{1}$, a contradiction. Thus $M$ and
$M^{d}=\left\{  f^{d}:f\in M\right\}  $ are disjoint.

The operation $f\ast g=f\left(  g,e_{2}\right)  $ is associative in any clone,
and $\left(  M;\ast\right)  $ is a semigroup in virtue of
\ref{lemma minimono ii}. The first projection is an identity element for
$\ast$, hence $\left(  M;\ast\right)  $ is a monoid. If $N$ is a submonoid of
$M$, then $N\cup N^{d} $ is closed under binary compositions. In a minimal
clone such a set must be either $\left\{  e_{1},e_{2}\right\}  $ or the whole
binary part of the clone. This fact together with the disjointness of $M$ and
$M^{d}$ shows that $Clo^{\left(  2\right)  }\left(  \mathbb{A}\right)  =M\cup
M^{d},$ and the only submonoids of $M$ are $\left\{  e_{1}\right\}  $ and $M$
itself. Such a monoid is called minimal, and it was shown in Claim 3.11 of
\cite{KSz comm} that every minimal monoid is isomorphic to a two-element
semilattice or a cyclic group of prime order.

Suppose first that $\left(  M;\ast\right)  \cong\left(  \left\{  0,1\right\}
;\vee\right)  $ with $f_{0}$ and $f_{1}$ corresponding to $0$ and $1$ at this
isomorphism. Then there are only four binary operations in $Clo\left(
\mathbb{A}\right)  $, namely $e_{1}=f_{0},e_{2}=f_{0}^{d},f_{1},f_{1}^{d}$ and
we can suppose (after passing to the dual of $\mathbb{A}$ if necessary) that
$f_{1}\left(  x,y\right)  =xy$, the basic operation in $\mathbb{A}$. By the
above isomorphism we have $f_{1}=f_{1\vee1}=f_{1}\ast f_{1}=f_{1}\left(
f_{1},e_{2}\right)  $, and this means that $xy\approx\left(  xy\right)  y$
holds in $\mathbb{A}$. Writing out \ref{lemma minimono i} with $f=f_{1}%
,g=f_{1},h=f_{0}$ and $f=f_{1},g=f_{0},h=f_{1}$ we get $f_{1}\left(
f_{1},f_{0}\right)  =f_{1}$ and $f_{1}\left(  f_{0},f_{1}\right)  =f_{0}$
implying that $\mathbb{A}$ satisfies the identities $\left(  xy\right)
x\approx xy$ and $x\left(  xy\right)  \approx x$. Similarly we obtain
$f_{1}\left(  f_{0},f_{1}^{d}\right)  =f_{1}\left(  f_{0},e_{2}\right)  $ and
$f_{1}\left(  f_{1},f_{1}^{d}\right)  =f_{1}\left(  f_{1},e_{2}\right)  $ as
special cases of \ref{lemma minimono ii}, and they translate to the identities
$x\left(  yx\right)  \approx xy$ and $\left(  xy\right)  \left(  yx\right)
\approx\left(  xy\right)  y$. All the identities in Proposition \ref{prop D}
are established, therefore $\mathbb{A\in}\mathcal{D}$ follows.

Now let us suppose that $\left(  M;\ast\right)  \cong\left(  \mathbb{Z}%
\hspace{0cm}_{p};+\right)  $ with $f_{i}\in M$ corresponding to $i\in
\mathbb{Z}\hspace{0cm}_{p}$. We have $f_{0}=e_{1}$ and we can suppose (after
dualizing if necessary) that $f_{i}\left(  x,y\right)  =xy$ for some
$i\in\left\{  1,\ldots,p-1\right\}  $. Since the automorphism group of
$\mathbb{Z}\hspace{0cm}_{p}$ acts transitively on $\left\{  1,\ldots
,p-1\right\}  $, we can suppose without loss of generality that $f_{1}\left(
x,y\right)  =xy$. Then $f_{i+1}=f_{1}\ast f_{i}=f_{1}\left(  f_{i}%
,e_{2}\right)  $, thus $f_{i+1}\left(  x,y\right)  =f_{i}\left(  x,y\right)
\cdot\nolinebreak y$, therefore $f_{i}\left(  x,y\right)  =xy^{i}$ and the
binary part of $Clo\left(  \mathbb{A}\right)  $ consists of the $2p$
operations $f_{i},f_{i}^{d}$ $\left(  i=0,1,\ldots,p-1\right)  $. Similarly to
the previous case, $\mathbb{F}_{2}\left(  V\left(  \mathbb{A}\right)  \right)
$ can be determined: \ref{lemma minimono i} implies $f_{i}\cdot f_{j}%
=f_{1}\left(  f_{i},f_{j}\right)  =f_{i}$, and \ref{lemma minimono ii} implies
$f_{i}\cdot f_{j}^{d}=f_{1}\left(  f_{i},f_{j}^{d}\right)  =f_{1}\left(
f_{i},e_{2}\right)  =f_{i+1}$; dualizing these we get $f_{i}^{d}\cdot
f_{j}^{d}=f_{i}^{d}$ and $f_{i}^{d}\cdot f_{j}=f_{i+1}^{d}$. It is easy to
check that $\mathbb{F}_{2}\left(  V\left(  \mathbb{A}\right)  \right)  $ is a
$p$-cyclic groupoid with a nontrivial clone (actually it is isomorphic to
$\mathbb{F}_{2}\left(  \mathcal{C}_{p}\right)  $), hence $V\left(
\mathbb{A}\right)  =\mathcal{C}_{p}$ by Proposition
\ref{prop p-cyclic, rect. band}.
\end{proof}

\begin{theorem}
\label{thm 1-2}Let $\mathcal{V}$ be a subvariety of $\mathcal{A}$ having a
minimal clone. Then $\mathcal{V}$ or its dual is a subvariety of
$\mathcal{B},\mathcal{C}_{p},\mathcal{D}$ or the variety of rectangular bands.
\end{theorem}

\begin{proof}
Let $\mathcal{W}$ be the intersection of $\mathcal{V}$ and the variety of
semigroups. Then $\mathcal{W}$ has a minimal or trivial clone, therefore it is
a subvariety of the variety of left zero semigroups, right zero semigroups,
rectangular bands, left regular bands or right regular bands (cf. \cite{PPP
prep},\cite{SzM minfcs}). We treat these five cases separately.

\begin{caselist}
\item \label{thm 1-2 Case 1}If $\mathcal{W}$ is the variety of left zero
semigroups, then Lemma \ref{lemma v-w} shows that $\mathcal{V}$ satisfies
$t_{1}x\approx t_{1}t$ for arbitrary terms $t_{1},t$ if $x$ is the first
variable of $t$. Specializing to $t=t_{1}$ we have that $\mathcal{V}\models
tx\approx tt\approx t$, i.e. a $\mathcal{V}$-term does not change if we
multiply it by its first variable from the right. Using these observations it
is easy to check that $M=\left\{  x,xy,xy^{2},xy^{3},\ldots\right\}  $
satisfies the conditions of Lemma \ref{lemma minimono} for any $\mathbb{A\in
}\mathcal{V}$ with a nontrivial clone (especially also for $\mathbb{F}%
_{\aleph_{0}}\left(  \mathcal{V}\right)  $), and hence $\mathcal{V}%
\subseteq\mathcal{D}$ or $\mathcal{V}=\mathcal{C}_{p}$ for some prime $p$.
(Note that $\mathcal{V}$ satisfies $x\left(  yz\right)  \approx xy$, therefore
Lemma 4.3 of \cite{W wa} could be used as well.)

\item \label{thm 1-2 Case 2}If $\mathcal{W}$ is the variety of right zero
semigroups, then similarly to the previous case we have the identities
$t_{1}x\approx t_{1}t$ and $tx\approx t$ in $\mathcal{V}$, where $x$ is the
last variable of $t$. Now we can apply Lemma \ref{lemma minimono} with
$\mathbb{A=F}_{\aleph_{0}}\left(  \mathcal{V}\right)  $ and $M=\left\{
x,\overleftarrow{xyx},\overleftarrow{xyxyx},\overleftarrow{xyxyxyx}%
,\ldots\right\}  $ to show that $\mathcal{V}\subseteq\mathcal{D}$ or
$\mathcal{V}=\mathcal{C}_{p}$ for some prime $p$, provided $\overleftarrow
{xyx}$ is nontrivial in $\mathbb{F}_{\aleph_{0}}\left(  \mathcal{V}\right)  $.
If $\left(  xy\right)  x$ is a projection in $\mathbb{F}_{\aleph_{0}}\left(
\mathcal{V}\right)  $, then $\mathcal{V}\models\left(  xy\right)  x\approx x$
or $\mathcal{V}\models\left(  xy\right)  x\approx y$. The latter is
impossible, since $x\left(  \left(  xy\right)  x\right)  \approx xx\approx x$
holds in $\mathcal{V}.$ Now we can write up the multiplication table of
$\mathbb{F}_{2}\left(  \mathcal{V}\right)  $.
\[%
\begin{tabular}
[c]{l|llll}%
$\cdot$ & $x$ & $y$ & $xy$ & $yx$\\\hline
$x$ & $x$ & $xy$ & $xy$ & $x$\\
$y$ & $yx$ & $y$ & $y$ & $yx$\\
$xy$ & $x$ & $xy$ & $xy$ & $x$\\
$yx$ & $yx$ & $y$ & $y$ & $yx$%
\end{tabular}
\]

This is a semigroup in $\mathcal{V}$, but it is not a right zero semigroup,
contradicting that $\mathcal{W}$ is the variety of right zero semigroups.
(Actually this groupoid is isomorphic to the two-generated free rectangular
band, hence Proposition \ref{prop p-cyclic, rect. band} could be applied as well.)

\item \label{thm 1-2 Case 3}If $\mathcal{W}$ is the variety of rectangular
bands, then $\mathcal{V}=\mathcal{W}$ by Proposition
\ref{prop p-cyclic, rect. band}.

\item \label{thm 1-2 Case 4}Suppose now that $\mathcal{W}$ is a variety of
left regular bands. Then $\mathcal{W}\models t_{1}\approx t_{2}$ if $t_{1}$
and $t_{2}$ are binary terms such that both $x$ and $y$ appear in both terms,
and they have the same first variable. Lemma \ref{lemma v-w} implies that
$tt_{1}\approx tt_{2}$ holds in $\mathcal{V}$ for every term $t$, if $t_{1}$
and $t_{2}$ satisfy the above conditions. This allows us to perform the
following computations in $\mathcal{V}$ with $g\left(  x,y\right)  =x\left(
xy\right)  $.%
\begin{align*}
g\left(  x,g\left(  x,y\right)  \right)   &  \approx x\left(  x\left(
x\left(  xy\right)  \right)  \right)  \approx x\left(  xy\right)  \approx
g\left(  x,y\right) \\
g\left(  x,g\left(  y,x\right)  \right)   &  \approx x\left(  x\left(
y\left(  yx\right)  \right)  \right)  \approx x\left(  xy\right)  \approx
g\left(  x,y\right) \\
g\left(  g\left(  x,y\right)  ,x\right)   &  \approx\left(  x\left(
xy\right)  \right)  \left(  \left(  x\left(  xy\right)  \right)  x\right)
\approx\left(  x\left(  xy\right)  \right)  \left(  x\left(  xy\right)
\right)  \approx g\left(  x,y\right) \\
g\left(  g\left(  x,y\right)  ,y\right)   &  \approx\left(  x\left(
xy\right)  \right)  \left(  \left(  x\left(  xy\right)  \right)  y\right)
\approx\left(  x\left(  xy\right)  \right)  \left(  x\left(  xy\right)
\right)  \approx g\left(  x,y\right) \\
g\left(  g\left(  x,y\right)  ,g\left(  y,x\right)  \right)   &
\approx\left(  x\left(  xy\right)  \right)  \left(  \left(  x\left(
xy\right)  \right)  \left(  y\left(  yx\right)  \right)  \right)
\approx\left(  x\left(  xy\right)  \right)  \left(  x\left(  xy\right)
\right)  \approx g\left(  x,y\right)
\end{align*}

These identities show that the subclone of $Clo\left(  \mathcal{V}\right)  $
generated by $g$ contains at most four binary operations, namely $g,g^{d}$ and
the two projections. If $g$ is nontrivial, then the minimality of the clone
implies that $g\left(  x,y\right)  =xy$ or $g\left(  y,x\right)  =yx.$ In the
first case the above identities are just the axioms of $\mathcal{B}$, and in
the second case they show that $\mathcal{V}\subseteq\mathcal{B}^{d}$. If $g$
is trivial, then $x\left(  xy\right)  \approx x$ holds in $\mathcal{V} $
(since $x\left(  xy\right)  \approx y$ is clearly impossible), and hence also
in $\mathcal{W}$. Since $\mathcal{W}$ is a variety of bands, $\mathcal{W}%
\models x\left(  xy\right)  \approx xy$, and therefore it is the variety of
left zero semigroups, and we have Case \ref{thm 1-2 Case 1}.

\item \label{thm 1-2 Case 5}Finally, let $\mathcal{W}$ be a variety of right
regular bands. Now $\mathcal{V}\models tt_{1}\approx tt_{2}$ whenever the last
variable of the binary terms $t_{1}$ and $t_{2}$ is the same, and the same
variables occur in them. Proceeding similarly to the previous case, we show
that $\left[  g\right]  ^{\left(  2\right)  }=\left\{  e_{1},e_{2}%
,g,g^{d}\right\}  $ for $g\left(  x,y\right)  =x\left(  yx\right)  $. This is
established by the following identities.%
\begin{align*}
g\left(  x,g\left(  x,y\right)  \right)   &  \approx x\left(  \left(  x\left(
yx\right)  \right)  x\right)  \approx x\left(  yx\right)  \approx g\left(
x,y\right)  \\
g\left(  x,g\left(  y,x\right)  \right)   &  \approx x\left(  \left(  y\left(
xy\right)  \right)  x\right)  \approx x\left(  yx\right)  \approx g\left(
x,y\right)  \\
g\left(  g\left(  x,y\right)  ,x\right)   &  \approx\left(  x\left(
yx\right)  \right)  \left(  x\left(  x\left(  yx\right)  \right)  \right)
\approx\left(  x\left(  yx\right)  \right)  \left(  x\left(  yx\right)
\right)  \approx g\left(  x,y\right)  \\
g\left(  g\left(  x,y\right)  ,y\right)   &  \approx\left(  x\left(
yx\right)  \right)  \left(  y\left(  x\left(  yx\right)  \right)  \right)
\approx\left(  x\left(  yx\right)  \right)  \left(  x\left(  yx\right)
\right)  \approx g\left(  x,y\right)  \\
g\left(  g\left(  x,y\right)  ,g\left(  y,x\right)  \right)   &
\approx\left(  x\left(  yx\right)  \right)  \left(  \left(  y\left(
xy\right)  \right)  \left(  x\left(  yx\right)  \right)  \right)
\approx\left(  x\left(  yx\right)  \right)  \left(  x\left(  yx\right)
\right)  \approx g\left(  x,y\right)
\end{align*}

If $g$ is nontrivial, then we have $\mathcal{V\subseteq B}$ or
$\mathcal{V\subseteq B}^{d}$ just as in Case \ref{thm 1-2 Case 4}. If $g$ is
trivial, then it has to be a first projection, hence $x\left(  yx\right)
\approx x$ holds in $\mathcal{V}$. Right regular bands satisfy $x\left(
yx\right)  \approx yx$, hence $\mathcal{W}\models yx\approx x$, and we have
Case \ref{thm 1-2 Case 2}.\qedhere
\end{caselist}
\end{proof}

Now we are ready to prove the main result of this section, the
characterization of groupoids with a minimal clone, that are almost semigroups
in the `spectral' sense.

\begin{theorem}
\label{THM spec}For a groupoid $\mathbb{A}$ the following two conditions are
equivalent%
%TCIMACRO{\TeXButton{roman}{\renewcommand{\theenumi}{(\roman{enumi})}
%\renewcommand{\labelenumi}{\theenumi}}}%
%BeginExpansion
\renewcommand{\theenumi}{(\roman{enumi})}
\renewcommand{\labelenumi}{\theenumi}%
%EndExpansion

\begin{enumerate}
\item \label{THM spec i}$\mathbb{A}$ has a minimal clone, and $1<s_{\mathbb{A}%
}\left(  4\right)  <5$;

\item \label{THM spec ii}$\mathbb{A}$ is not a semigroup, and $\mathbb{A}$ or
its dual belongs to one of the varieties $\mathcal{B}\cap\mathcal{A}$,
$\mathcal{C}_{p}$, or $\mathcal{D}\cap\mathcal{A}$.
\end{enumerate}

\noindent If these conditions are fulfilled, then we have $s_{\mathbb{A}%
}\left(  n\right)  =2^{n-2}$ for $n\geq2$.
\end{theorem}

\begin{proof}
First we show that \ref{THM spec i} implies \ref{THM spec ii}. The
considerations preceding Lemma \ref{lemma general ass} show that if
$\mathbb{A}$ has a minimal clone, and $1<s_{\mathbb{A}}\left(  4\right)  <5,$
then either $\mathbb{A}$ or its dual satisfies $x\left(  y\left(  zu\right)
\right)  \approx x\left(  \left(  yz\right)  u\right)  $, i.e. $\mathbb{A}\in$
$\mathcal{A}$ or $\mathbb{A}\in\mathcal{A}^{d}$. Applying Theorem
\ref{thm 1-2}, we get that $\mathbb{A}$ or $\mathbb{A}^{d}$ belongs to
$\mathcal{B},\mathcal{C}_{p}$ or $\mathcal{D}$ (for some prime $p$). Thus we
have to consider varieties of the from $\mathcal{V}_{1}\cap\mathcal{V}_{2}$,
where $\mathcal{V}_{1}=\mathcal{A}$ or $\mathcal{V}_{1}=\mathcal{A}^{d}$, and
$\mathcal{V}_{2}\in\left\{  \mathcal{B},\mathcal{C}_{p},\mathcal{D}%
,\mathcal{B}^{d},\mathcal{C}_{p}^{d},\mathcal{D}^{d}:p\text{ is a
prime}\right\}  $, but up to duality we have only six cases, because we may
suppose that $\mathcal{V}_{2}=\mathcal{B},\mathcal{C}_{p}$ or $\mathcal{D}$.

We show that if $\mathbb{A}\in\mathcal{V}_{2}$, and $a,b$ are elements of
$\mathbb{A}$ such that $ax=bx$ holds for all $x\in\mathbb{A}$, then $a=b$.
Letting $x=a$ and $x=b$ we see that $\left\{  a,b\right\}  $ is a right zero
subsemigroup of $\mathbb{A}$. The identity $x\left(  yx\right)  \approx xy$
holds in $\mathcal{V}_{2}$ in all of the three cases, hence $a\left(
ba\right)  =ab$. Since $a$ and $b$ form a right zero semigroup we have
$a\left(  ba\right)  =a$ and $ab=b$, thus $a=b$ as claimed. We see that
$\mathcal{V}_{2}\cap\mathcal{A}^{d}$ is a variety of semigroups, because the
defining identity of $\mathcal{A}^{d}$ is $\left(  \left(  xy\right)
z\right)  u\approx\left(  x\left(  yz\right)  \right)  u$, and according to
the previous observation this implies that $\left(  xy\right)  z\approx
x\left(  yz\right)  $ holds in $\mathcal{V}_{2}$. Thus $\mathcal{V}%
_{1}=\mathcal{A}$, and we end up with the varieties of \ref{THM spec ii}.
(Note that $\mathcal{C}_{p}\models x\left(  y\left(  zu\right)  \right)
\approx xy\approx x\left(  \left(  yz\right)  u\right)  $, therefore
$\mathcal{C}_{p}\cap\mathcal{A=C}_{p}$.)

Now suppose that $\mathbb{A}$ (or its dual) belongs to one of the varieties
mentioned in \ref{THM spec ii}, and $\mathbb{A}$ is not a semigroup. The clone
of $\mathcal{B}$, $\mathcal{C}_{p}$ and $\mathcal{D}$ is minimal, thus the
clone of $\mathbb{A}$ is minimal, too (note that $\mathbb{A}$ has a nontrivial
clone, because it is not a semigroup). The other assertion of \ref{THM spec i}
will follow at once, if we prove that $s_{\mathbb{A}}\left(  n\right)
=2^{n-2}$. We will do this in two steps: first we show that $\mathbb{A}%
\in\mathcal{A}$ implies $s_{\mathbb{A}}\left(  n\right)  \leq2^{n-2}$, and
then we prove that $s_{\mathbb{A}}\left(  n\right)  \geq2^{n-2}$ holds if we
suppose in addition that $\mathbb{A}\in\mathcal{B},\mathcal{C}_{p}$ or
$\mathcal{D}$.

Let $B$ and $B^{\prime}$ be bracketings of the product $x_{1}\cdot\ldots\cdot
x_{n}$. Lemma \ref{lemma general ass} implies that $\mathcal{A}\models
B\approx B^{\prime}$ if $\left\vert l\left(  B\right)  \right\vert =\left\vert
l\left(  B^{\prime}\right)  \right\vert $ and $\mathcal{A}\models l\left(
B\right)  \approx l\left(  B^{\prime}\right)  $. Applying Lemma
\ref{lemma general ass} again, we see that $\left\vert l\left(  B\right)
\right\vert =\left\vert l\left(  B^{\prime}\right)  \right\vert ,\left\vert
l^{2}\left(  B\right)  \right\vert =\left\vert l^{2}\left(  B^{\prime}\right)
\right\vert $ and $l^{2}\left(  B\right)  \approx l^{2}\left(  B^{\prime
}\right)  $ is sufficient for $B\approx B^{\prime}$. Proceeding this way we
arrive at left factors of size $1$ (i.e. the single variable $x_{1}$) finally,
and we see that if $\left\vert l^{i}\left(  B\right)  \right\vert =\left\vert
l^{i}\left(  B^{\prime}\right)  \right\vert $ for all $i$ (where it makes
sense), then $B\approx B^{\prime}$ holds in $\mathcal{A}$. Clearly, the
numbers $\left\vert l^{i}\left(  B\right)  \right\vert $ (and $\left\vert
l^{i}\left(  B^{\prime}\right)  \right\vert $) are strictly decreasing in $i$,
therefore it is sufficient if the sets $\left\{  \left\vert l^{i}\left(
B\right)  \right\vert :i=1,2,\ldots\right\}  $ and $\left\{  \left\vert
l^{i}\left(  B^{\prime}\right)  \right\vert :i=1,2,\ldots\right\}  $ coincide.
They are subsets of $\left\{  1,2,\ldots,n-1\right\}  $, containing $1$, hence
there are $2^{n-2}$ many choices for these sets. This shows that
$s_{\mathbb{A}}\left(  n\right)  \leq2^{n-2}$ for any $\mathbb{A}%
\in\mathcal{A}$.

Now let $\mathbb{A}\in\mathcal{A}\cap\mathcal{V}_{2}$, where $\mathcal{V}%
_{2}\in\left\{  \mathcal{B},\mathcal{C}_{p},\mathcal{D}:p\text{ is a
prime}\right\}  $, and let $B$ and $B^{\prime}$ be bracketings as before.
Suppose that $\mathbb{A}\models B\approx B^{\prime}$, but $\left\{  \left\vert
l^{i}\left(  B\right)  \right\vert :i=1,2,\ldots\right\}  \neq\left\{
\left\vert l^{i}\left(  B^{\prime}\right)  \right\vert :i=1,2,\ldots\right\}
$, and let $i$ be the smallest value where $\left\vert l^{i}\left(  B\right)
\right\vert $ and $\left\vert l^{i}\left(  B^{\prime}\right)  \right\vert $
are different. Applying Lemma \ref{lemma general ass} and the observation made
in the second paragraph of this proof (a certain right cancellation property)
we can delete the right factors in the identity $B\approx B^{\prime}$, if they
have the same size. Doing this $i-1$ times we arrive at bracketings whose left
factors have different size, thus we may suppose that $i=1$, and we can also
suppose that $\left\vert l^{1}\left(  B\right)  \right\vert <\left\vert
l^{1}\left(  B^{\prime}\right)  \right\vert $. Let us substitute $x$ for the
first $\left\vert l^{1}\left(  B\right)  \right\vert $ variables, $y$ for the
next $\left\vert l^{1}\left(  B^{\prime}\right)  \right\vert -\left\vert
l^{1}\left(  B\right)  \right\vert $ variables, and $z$ for the rest. Then $B$
becomes $\left(  x\cdots x\right)  \left(  y\cdots yz\cdots z\right)  $ (with
some bracketing of the two products), and $B^{\prime}$ has the form $\left(
x\cdots xy\cdots y\right)  \left(  z\cdots z\right)  $. Thus $\mathbb{A}$
satisfies an identity of the form $\left(  x\cdots x\right)  \left(  y\cdots
yz\cdots z\right)  \approx\left(  x\cdots xy\cdots y\right)  \left(  z\cdots
z\right)  $ (with the same number of $x$, $y$ and $z$ on the two sides).

In $\mathcal{B}$ this identity reduces to $x\left(  yz\right)  \approx\left(
xy\right)  z$, showing that if $s_{\mathbb{A}}\left(  n\right)  <2^{n-2}$ for
some $n$, then $\mathbb{A}$ is a semigroup. If $\mathcal{V}_{2}=\mathcal{C}%
_{p}$ or $\mathcal{D}$, then let us put $y=x$, then we have $\mathbb{A}%
\models\left(  x\cdots x\right)  \left(  x\cdots xz\cdots z\right)
\approx\left(  x\cdots xx\cdots x\right)  \left(  z\cdots z\right)  $. The
right hand side is clearly $xz$, and on the left hand side the bracketing of
the factor $\left(  x\cdots xz\cdots z\right)  $ is irrelevant according to
Lemma \ref{lemma general ass}. Thus $\mathbb{A}\models x\left(  xz\right)
\approx xz$, and since $x\left(  xz\right)  \approx x$ holds in $\mathcal{C}%
_{p}$ and $\mathcal{D}$ we see that $\mathbb{A}$ is a left zero semigroup. We
have proved that the associative spectrum of a groupoid in any one of the
varieties mentioned in \ref{THM spec ii} is either $\left(  1,1,1,1,\ldots
\right)  $ or $\left(  1,2,4,8,\ldots\right)  $, and this completes the proof
of the theorem.
\end{proof}

\begin{remark}
Each of the varieties $\mathcal{B}\cap\mathcal{A}$, $\mathcal{C}_{p}$ and
$\mathcal{D}\cap\mathcal{A}$ contain groupoids with a nonassociative
operation. For $\mathcal{C}_{p}$ it is clear, because the only $p$-cyclic
groupoids that are semigroups are the left zero semigroups. The two-generated
free algebra of $\mathcal{D}$ is not a semigroup, and satisfies $x\left(
y\left(  zu\right)  \right)  \approx x\left(  \left(  yz\right)  u\right)  $,
hence belongs to $\mathcal{D}\cap\mathcal{A}$. (See the multiplication table
in the proof of Proposition \ref{prop D}.) Let us now construct some
nonassociative algebras in $\mathcal{B}\cap\mathcal{A}$.

Let $\mathbb{S}=\left(  S;\vee\right)  $ be a semilattice, and let $C$ be the
set of finite chains in $\mathbb{S}$. We define a multiplication in $C$ by the
following formula (note that if $b_{l}\leq a_{k}$, then the right hand side is
the same as the first factor on the left hand side).%
\[
\left(  a_{1}<a_{2}<\cdots<a_{k}\right)  \cdot\left(  b_{1}<b_{2}<\cdots
<b_{l}\right)  =\left(  a_{1}<a_{2}<\cdots<a_{k}\leq a_{k}\vee b_{l}\right)
\]
For $\mathbf{a=}\left(  a_{1}<a_{2}<\cdots<a_{k}\right)  ,\mathbf{b=}\left(
b_{1}<b_{2}<\cdots<b_{l}\right)  $ and $\mathbf{c=}\left(  c_{1}<c_{2}%
<\cdots<c_{m}\right)  $ we have $\left(  \mathbf{a}\cdot\mathbf{b}\right)
\cdot\mathbf{c=}\left(  a_{1}<a_{2}<\cdots<a_{k}\leq a_{k}\vee b_{l}\leq
a_{k}\vee b_{l}\vee c_{m}\right)  $ and $\mathbf{a}\cdot\left(  \mathbf{b}%
\cdot\mathbf{c}\right)  =\left(  a_{1}<a_{2}<\cdots<a_{k}\leq a_{k}\vee
b_{l}\vee c_{m}\right)  $. Since the top element of both chains is $a_{k}\vee
b_{l}\vee c_{m}$, right multiplication by $\left(  \mathbf{a}\cdot
\mathbf{b}\right)  \cdot\mathbf{c}$ is the same as right multiplication by
$\mathbf{a}\cdot\left(  \mathbf{b}\cdot\mathbf{c}\right)  $, hence
$\mathbb{C}=\left(  C;\cdot\right)  $ satisfies $x\left(  y\left(  zu\right)
\right)  \approx x\left(  \left(  yz\right)  u\right)  $. It is not hard to
check, that the defining identities of $\mathcal{B}$ also hold in $\mathbb{C}%
$, hence $\mathbb{C}\in\mathcal{B}\cap\mathcal{A}$. If the height of
$\mathbb{S}$ is at least three, i.e. there is a chain of length three, then
$\mathbb{C}$ is not a semigroup. Indeed, if $a<b<c$, then $\left(  a\cdot
b\right)  \cdot c=\left(  a<b<c\right)  \neq\left(  a<c\right)  =a\cdot\left(
b\cdot c\right)  $.
\end{remark}

\begin{remark}
The variety $\mathcal{D}\cap\mathcal{A}$ was defined by an infinite set of
identities, but it has a finite basis, namely $xx\approx x,x\left(  yz\right)
\approx xy,\left(  xy\right)  y\approx xy$. Indeed, it is quite
straightforward to check that any algebra satisfying these identities belongs
to $\mathcal{D}\cap\mathcal{A}$. Conversely, if $\mathbb{A\in}\mathcal{D}%
\cap\mathcal{A}$, then $\mathbb{A}\models x\left(  yz\right)  \approx x\left(
\left(  yy\right)  z\right)  \approx x\left(  y\left(  yz\right)  \right)
\approx xy$, and $\mathbb{A}$ also satisfies $xx\approx x$ and $\left(
xy\right)  y\approx xy$ as they are among the defining identities of
$\mathcal{D}$. This latter axiomatization of $\mathcal{D}\cap\mathcal{A}$
resembles to the definition of $p$-cyclic groupoids. It is an interesting fact
that every groupoid that has a minimal clone and satisfies $x\left(
yz\right)  \approx xy$ belongs to one of the varieties $\mathcal{C}_{p}$ or
$\mathcal{D}\cap\mathcal{A}$ (cf. Lemma 4.3 of \cite{W wa}).
\end{remark}

\section{Sz\'{a}sz-H\'{a}jek groupoids with a minimal clone\label{Section SH}}

Another way to measure associativity is to count the number of
\emph{nonassociative triples} in the groupoid; this number (or cardinal, in
the infinite case) is called the \emph{index of nonassociativity}, and is
denoted by $ns$. Formally, we have $ns\left(  \mathbb{A}\right)  =\left\vert
\left\{  \left(  a,b,c\right)  \in A^{3}:\left(  ab\right)  c\neq a\left(
bc\right)  \right\}  \right\vert $. This notion has been studied by several
authors \cite{Climescu1, Climescu2, DrapalKepka, KT1 ns, Szasz}. Clearly
$\mathbb{A}$ is a semigroup iff $ns\left(  \mathbb{A}\right)  =0$, and it is
natural to say that the multiplication of $\mathbb{A}$ is almost associative,
if $ns\left(  \mathbb{A}\right)  =1$. Such groupoids are called
Sz\'{a}sz-H\'{a}jek groupoids (SH-groupoids for short). SH-groupoids were
investigated in \cite{Hajek1, Hajek2} and \cite{KT3 SH aaa, KT4 aba, KT5 aab,
KT6 abc} in much detail. Following the terminology of these papers, we say
that an SH-groupoid is of type $\left(  a,b,c\right)  $, if its only
nonassociative triple is $\left(  a,b,c\right)  \in A^{3}$ and $a\neq b\neq
c\neq a$. Types $\left(  a,a,a\right)  ,\left(  a,b,a\right)  ,\left(
a,a,b\right)  $ and $\left(  a,b,b\right)  $ are defined analogously. (Note
that by saying e.g. that $\mathbb{A}$ is an SH-groupoid of type $\left(
a,b,c\right)  $ we mean not only that the components of the unique
nonassociative triple are pairwise distinct, but implicitly we assume that
these components are denoted by $a$, $b$ and $c$ respectively.) Let us recall
a result from \cite{KT3 SH aaa} (Proposition 1.2(i)).

\begin{proposition}
\label{prop xy=a}If $\mathbb{A}$ is an SH-groupoid, and $\left(  a,b,c\right)
$ is the unique nonassociative triple, then $xy=a$ $(xy=b,xy=c)$ implies $x=a
$ $\left(  x=b,x=c\right)  $ or $y=a$ $\left(  y=b,y=c\right)  $ for all
$x,y\in A$.
\end{proposition}

\begin{proof}
Suppose that $xy=a$, but $x\neq a\neq y$. Since $x\neq a$, we have $\left(
x,y,bc\right)  \neq\left(  a,b,c\right)  $, hence $\left(  x,y,bc\right)  $ is
an associative triple: $\left(  xy\right)  \left(  bc\right)  =x\left(
y\left(  bc\right)  \right)  $. Now $y\neq a$ implies that $\left(
y,b,c\right)  \neq\left(  a,b,c\right)  $, so $x\left(  y\left(  bc\right)
\right)  =x\left(  \left(  yb\right)  c\right)  $. Similarly $x\left(  \left(
yb\right)  c\right)  =\left(  x\left(  yb\right)  \right)  c=\left(  \left(
xy\right)  b\right)  c$, because $x\neq a.$ We have obtained that $\left(
xy\right)  \left(  bc\right)  =\left(  \left(  xy\right)  b\right)  c$, thus
$\left(  xy,b,c\right)  =\left(  a,b,c\right)  $ is an associative triple,
which is a contradiction. The other two assertions can be proved similarly.
\end{proof}

Clearly, a subgroupoid of an SH-groupoid $\mathbb{A}$ with nonassociative
triple $\left(  a,b,c\right)  $ is an SH-groupoid or a semigroup, depending on
whether it contains $a$, $b$ and $c$ or not. Specially, $\mathbb{A}$ is
generated by $\left\{  a,b,c\right\}  $ iff all proper subgroupoids of
$\mathbb{A}$ are semigroups. Such a groupoid is called a \emph{minimal
SH-groupoid}. In \cite{KT3 SH aaa, KT4 aba, KT5 aab, KT6 abc} the project of
characterizing minimal SH-groupoids was begun, but completed only for the type
$\left(  a,a,a\right)  $. In Theorem \ref{THM SH} we prove that SH-groupoids
having a minimal clone belong to the varieties $\mathcal{B}$ or $\mathcal{B}%
^{d}$, and in Theorem \ref{thm minimal SH in B} we give a complete list of
minimal SH-groupoids with a minimal clone up to isomorphism. We need the
following lemma before we state and prove the main result.

\begin{lemma}
If an SH-groupoid has a minimal clone, then it is of type $\left(
a,b,c\right)  $.
\end{lemma}

\begin{proof}
Let $\mathbb{A}$ be an SH-groupoid with a minimal clone. Then $\mathbb{A}$ is
idempotent, hence it cannot be of type $\left(  a,a,a\right)  $. If it is of
type $\left(  a,b,a\right)  $, then the subgroupoid generated by $a$ and $b$
is a minimal SH-groupoid of type $\left(  a,b,a\right)  $ with a minimal
clone. The description of minimal SH-groupoids of type $\left(  a,b,a\right)
$ given in \cite{KT4 aba} is not complete, but it covers the idempotent case
(subtypes $\left(  \alpha\right)  $ and $\left(  \beta\right)  $). There are
four idempotent minimal SH-groupoids of type $\left(  a,b,a\right)  $ up to
isomorphism: the following two groupoids and their duals (the second groupoid
is a factor of the first one).%

\[%
\begin{tabular}
[c]{l|llll}%
$\cdot$ & $a$ & $b$ & $d$ & $e$\\\hline
$a$ & $a$ & $a$ & $e$ & $e$\\
$b$ & $d$ & $b$ & $d$ & $d$\\
$d$ & $d$ & $d$ & $d$ & $d$\\
$e$ & $e$ & $e$ & $e$ & $e$%
\end{tabular}
\qquad%
\begin{tabular}
[c]{llll}%
$\cdot$ & \multicolumn{1}{|l}{$a$} & $b$ & $d$\\\hline
$a$ & \multicolumn{1}{|l}{$a$} & $a$ & $d$\\
$b$ & \multicolumn{1}{|l}{$d$} & $b$ & $d$\\
$d$ & \multicolumn{1}{|l}{$d$} & $d$ & $d$\\
&  &  &
\end{tabular}
\]

In both cases the operation $g\left(  x,y\right)  =x\left(  yx\right)  $ is
nontrivial, and preserves the equivalence relation corresponding to the
partition whose only nontrivial block is $\{b,d\}$, but the basic operation
$f\left(  x,y\right)  =xy$ does not preserve this relation. This shows that
$f\notin\left[  g\right]  $, hence the clone is not minimal.

Suppose now that $\mathbb{A}$ is of type $\left(  a,a,b\right)  $. From the
computations in \cite{KT5 aab} it follows that $d=ba=b$ (combine Lemmas 1.5,
1.6, 2.4 and 2.19), therefore the subgroupoid generated by $a$ and $b$ is a
minimal SH-groupoid of type $\left(  a,a,b\right)  $ and of subtype $\left(
\varepsilon\right)  $. Up to isomorphism there is only one such groupoid,
namely the following one.
\[%
\begin{tabular}
[c]{l|llll}%
$\cdot$ & $a$ & $b$ & $c$ & $e$\\\hline
$a$ & $a$ & $c$ & $e$ & $e$\\
$b$ & $b$ & $b$ & $b$ & $b$\\
$c$ & $c$ & $c$ & $c$ & $c$\\
$e$ & $e$ & $e$ & $e$ & $e$%
\end{tabular}
\]

The clone of this groupoid is not minimal, because $x\left(  xy\right)  $ is a
nontrivial operation preserving the set $\{a,b,e\}$, while the basic operation
$xy$ does not preserve this set.

Dually, the type $\left(  a,b,b\right)  $ is not possible either, thus we can
conclude that an SH-groupoid with a minimal clone has to be of type $\left(
a,b,c\right)  $.
\end{proof}

\begin{theorem}
\label{THM SH}For a Sz\'{a}sz-H\'{a}jek groupoid $\mathbb{A}$ the following
two conditions are equivalent.%
%TCIMACRO{\TeXButton{roman}{\renewcommand{\theenumi}{(\roman{enumi})}
%\renewcommand{\labelenumi}{\theenumi}}}%
%BeginExpansion
\renewcommand{\theenumi}{(\roman{enumi})}
\renewcommand{\labelenumi}{\theenumi}%
%EndExpansion

\begin{enumerate}
\item \label{THM SH i}$\mathbb{A}$ has a minimal clone;

\item \label{THM SH ii}$\mathbb{A}$ or its dual belongs to the variety
$\mathcal{B}$.
\end{enumerate}
\end{theorem}

\begin{proof}
It is clear that \ref{THM SH ii} implies \ref{THM SH i}, since $\mathcal{B}$
has a minimal clone. For the other direction let us suppose that $\mathbb{A}$
is an SH-groupoid with a minimal clone. As we have seen in the previous lemma,
$\mathbb{A}$ is of type $\left(  a,b,c\right)  $. Therefore $\left(
x,x,y\right)  $ is an associative triple for all $x,y\in A$, hence
$\mathbb{A}\models x\left(  xy\right)  \approx xy$ by idempotence. Similarly,
we obtain $\mathbb{A}\models\left(  xy\right)  y\approx xy$ and $\mathbb{A}%
\models x\left(  yx\right)  \approx\left(  xy\right)  x$. Proposition
\ref{prop xy=a} shows that $\left(  x,y,xy\right)  $ is an associative triple
for all $x,y\in A $, because $x=a,y=b,xy=c$ is impossible. Thus $\mathbb{A}%
\models x\left(  y\left(  xy\right)  \right)  \approx\left(  xy\right)
\left(  xy\right)  \approx xy$. By another application of Proposition
\ref{prop xy=a} we can see that $\left(  xy,y,x\right)  \neq\left(
a,b,c\right)  $, so $\left(  xy\right)  \left(  yx\right)  \approx\left(
\left(  xy\right)  y\right)  x\approx\left(  xy\right)  x$ holds in
$\mathbb{A}$.

The identities derived so far are almost sufficient to fill out the
multiplication table of the two-generated free algebra in the variety
generated by $\mathbb{A}$ (see the table below). The only entries that are not
determined yet are $\left(  xyx\right)  \left(  yxy\right)  $ and $\left(
yxy\right)  \left(  xyx\right)  $. In order to compute these, let us observe
that $\left(  xyx,yx,y\right)  $ is always an associative triple, because
$yx=b $ and $y=c$ implies $x=b$ by Proposition \ref{prop xy=a}, but then
$xyx=bb=b\neq a$. Therefore $\mathbb{A}\models\left(  xyx\right)  \left(
yxy\right)  \approx\left(  \left(  xyx\right)  \left(  yx\right)  \right)
y\approx\left(  xyx\right)  y\approx xy$.%

\[%
\begin{tabular}
[c]{l|llllll}%
$\cdot$ & $x$ & $y$ & $xy$ & $yx$ & $xyx$ & $yxy$\\\hline
$x$ & $x$ & $xy$ & $xy$ & $xyx$ & $xyx$ & $xy$\\
$y$ & $yx$ & $y$ & $yxy$ & $yx$ & $yx$ & $yxy$\\
$xy$ & $xyx$ & $xy$ & $xy$ & $xyx$ & $xyx$ & $xy$\\
$yx$ & $yx$ & $yxy$ & $yxy$ & $yx$ & $yx$ & $yxy$\\
$xyx$ & $xyx$ & $xy$ & $xy$ & $xyx$ & $xyx$ & $xy$\\
$yxy$ & $yx$ & $yxy$ & $yxy$ & $yx$ & $yx$ & $yxy$%
\end{tabular}
\]

We see that the binary part of $Clo\left(  \mathbb{A}\right)  $ contains at
most six operations (some of the six elements in the table may coincide). In
\cite{LLPPP} we can find the complete description of minimal clones with at
most six binary operations, so we could finish the proof by simply examining
the list of clones given there.

Another way is to observe that for $g\left(  x,y\right)  =xyx$ the binary part
of $\left[  g\right]  $ is $\left\{  e_{1},e_{2},g,g^{d}\right\}  $. If $g$ is
a nontrivial operation, then $\left[  g\right]  =Clo\left(  \mathbb{A}\right)
$, hence $\mathbb{A}$ satisfies $xyx\approx xy$ or $xyx\approx yx$, and then
the defining identities of $\mathcal{B}$ or $\mathcal{B}^{d}$ can be read from
the above multiplication table. If $g$ is trivial, then $\mathbb{A}\models
xyx\approx x$, because $xyx\approx y$ would imply $xy\approx\left(
xyx\right)  y\approx yy\approx y$. In this case $\mathbb{F}_{2}\left(
V\left(  \mathbb{A}\right)  \right)  $ is a rectangular band (we get the same
multiplication table as in Case \ref{thm 1-2 Case 2} of the proof of Theorem
\ref{thm 1-2}), hence $\mathbb{A}$ is a rectangular band by Proposition
\ref{prop p-cyclic, rect. band}, contradicting that $\mathbb{A}$ is an SH-groupoid.
\end{proof}

Finally we describe minimal SH-groupoids in the varieties $\mathcal{B}$ and
$\mathcal{B}^{d}$ up to isomorphism.

\begin{theorem}
\label{thm minimal SH in B}Every minimal SH-groupoid having a minimal clone is
isomorphic or dually isomorphic to one of the groupoids $\mathbb{G}_{1}%
,\ldots,\mathbb{G}_{10}$ (see the multiplication tables in the proof).
\end{theorem}

\begin{proof}
Let $\mathbb{A}$ be a minimal SH-groupoid with a minimal clone. Then
$\mathbb{A}$ is of type $\left(  a,b,c\right)  $, and up to duality we may
suppose that $\mathbb{A}$ belongs to the variety $\mathcal{B}$. Following the
notation of \cite{KT6 abc} we set $d=ab,e=bc,f=a\left(  bc\right)  =ae$ and
$g=\left(  ab\right)  c=dc$. Some of these elements may coincide, but $a,b,c$
are pairwise distinct and $f\neq g$. Since $\mathbb{A}$ is idempotent, we have
$d=a$ or $e=c$ by Lemma 1.7 of \cite{KT6 abc}. If $d=a$, then $ba=b$ or $ba=a$
(Lemma 1.9 (iii)); if $e=c$, then $cb=b$ or $cb=c$ (Lemma 1.9 (iv)). Thus we
have four cases, and we will deal with them separately.

\begin{caselist}
\item $d=ab=a$ and $ba=b$

We have $g=dc=ac=c$ by Lemma 1.4 (ii) of \cite{KT6 abc}, and then $ca=c\left(
ca\right)  =\left(  ac\right)  \left(  ca\right)  =ac=c$ follows applying the
defining identities of $\mathcal{B}$. Some other products may be computed with
the help of these identities, for example $be=b\left(  bc\right)  =bc=e$ and
$eb=\left(  bc\right)  b=bc=e$. For others, we can use the fact that $\left(
a,b,c\right)  $ is the only nonassociative triple, e.g.: $cb=\left(
ca\right)  b=c\left(  ab\right)  =ca=c$, and $bf=b\left(  ae\right)  =\left(
ba\right)  e=be=e$.

We can fill out the multiplication table this way except for the entry $fc$.
Here we have two possibilities. If $f\neq a$ or $e\neq b$, then $\left(
f,e,c\right)  \neq\left(  a,b,c\right)  $, therefore $fc=\left(  fe\right)
c=f\left(  ec\right)  =fe=f$, and we get the groupoid $\mathbb{G}_{1}$. If
$f=a $ and $e=b$, then $fc=ac=c$, and we arrive at the groupoid $\mathbb{G}%
_{3}$. In both cases we have to consider the possibility that some of the
elements (denoted by different letters so far) coincide. This amounts to
forming factor groupoids, but only with respect to congruences where $f$ and
$g$ belong to different congruence classes (otherwise the factor groupoid
would be a semigroup). There is no such congruence on $\mathbb{G}_{3}$, while
$\mathbb{G}_{1}$ has exactly one nontrivial congruence not collapsing $f$ and
$g\left(  =c\right)  $; its classes are $\left\{  a\right\}  ,\left\{
b\right\}  ,\left\{  c\right\}  ,\left\{  e,f\right\}  $, and the
corresponding factor groupoid is $\mathbb{G}_{2}$.%
\[
\mathbb{G}_{1}:\
\begin{tabular}
[c]{l|lllll}%
$\cdot$ & $a$ & $b$ & $c$ & $e$ & $f$\\\hline
$a$ & $a$ & $a$ & $c$ & $f$ & $f$\\
$b$ & $b$ & $b$ & $e$ & $e$ & $e$\\
$c$ & $c$ & $c$ & $c$ & $c$ & $c$\\
$e$ & $e$ & $e$ & $e$ & $e$ & $e$\\
$f$ & $f$ & $f$ & $f$ & $f$ & $f$%
\end{tabular}
\ \qquad\mathbb{G}_{2}:\
\begin{tabular}
[c]{lllll}%
$\cdot$ & \multicolumn{1}{|l}{$a$} & $b$ & $c$ & $e$\\\hline
$a$ & \multicolumn{1}{|l}{$a$} & $a$ & $c$ & $e$\\
$b$ & \multicolumn{1}{|l}{$b$} & $b$ & $e$ & $e$\\
$c$ & \multicolumn{1}{|l}{$c$} & $c$ & $c$ & $c$\\
$e$ & \multicolumn{1}{|l}{$e$} & $e$ & $e$ & $e$\\
&  &  &  &
\end{tabular}
\]%
\[
\mathbb{G}_{3}:\
\begin{tabular}
[c]{l|lll}%
$\cdot$ & $a$ & $b$ & $c$\\\hline
$a$ & $a$ & $a$ & $c$\\
$b$ & $b$ & $b$ & $b$\\
$c$ & $c$ & $c$ & $c$%
\end{tabular}
\]

\item $d=ab=a$ and $ba=a$

Let us start again with the product $ca$. We claim that $\left(
a,b,ca\right)  $ is a nonassociative triple. Indeed, $\left(  ab\right)
\left(  ca\right)  =a\left(  ca\right)  =ac=\left(  ab\right)  c$, while
$a\left(  b\left(  ca\right)  \right)  =a\left(  \left(  bc\right)  a\right)
=a\left(  ea\right)  =ae=a\left(  bc\right)  $. Since the only nonassociative
triple is $\left(  a,b,c\right)  $, we can conclude that $ca=c$. Then
$cb=\left(  ca\right)  b=c\left(  ab\right)  =ca=c$, and the rest of the
multiplication table can be filled out without any difficulty (we will not
have to deal with a situation like that of $fc$ in the previous case). We get
the groupoid $\mathbb{G}_{4}$, and the only possible coincidence between the
six elements is $e=f$; this yields $\mathbb{G}_{5}$.%
\[
\mathbb{G}_{4}:\
\begin{tabular}
[c]{l|llllll}%
$\cdot$ & $a$ & $b$ & $c$ & $e$ & $f$ & $g$\\\hline
$a$ & $a$ & $a$ & $g$ & $f$ & $f$ & $g$\\
$b$ & $a$ & $b$ & $e$ & $e$ & $f$ & $g$\\
$c$ & $c$ & $c$ & $c$ & $c$ & $c$ & $c$\\
$e$ & $e$ & $e$ & $e$ & $e$ & $e$ & $e$\\
$f$ & $f$ & $f$ & $f$ & $f$ & $f$ & $f$\\
$g$ & $g$ & $g$ & $g$ & $g$ & $g$ & $g$%
\end{tabular}
\ \qquad\mathbb{G}_{5}:\
\begin{tabular}
[c]{llllll}%
$\cdot$ & \multicolumn{1}{|l}{$a$} & $b$ & $c$ & $e$ & $g$\\\hline
$a$ & \multicolumn{1}{|l}{$a$} & $a$ & $g$ & $e$ & $g$\\
$b$ & \multicolumn{1}{|l}{$a$} & $b$ & $e$ & $e$ & $g$\\
$c$ & \multicolumn{1}{|l}{$c$} & $c$ & $c$ & $c$ & $c$\\
$e$ & \multicolumn{1}{|l}{$e$} & $e$ & $e$ & $e$ & $e$\\
$g$ & \multicolumn{1}{|l}{$g$} & $g$ & $g$ & $g$ & $g$\\
&  &  &  &  &
\end{tabular}
\]

\item $e=bc=c$ and $cb=b$

As the following computation shows, this case is not possible, because the
identities of $\mathcal{B}$ imply that $\mathbb{A}$ is a semigroup. (We have
indicated where we used the axioms of $\mathcal{B}$ and the
Sz\'{a}sz-H\'{a}jek property.)
\begin{multline*}
g=dc\overset{\mathcal{B}}{=}\left(  dc\right)  c=\left(  dc\right)
e\overset{SH}{=}d\left(  ce\right)  =d\left(  c\left(  bc\right)  \right) \\
\overset{\mathcal{B}}{=}d\left(  cb\right)  =\left(  ab\right)  \left(
cb\right)  \overset{SH}{=}a\left(  b\left(  cb\right)  \right)  \overset
{\mathcal{B}}{\mathcal{=}}a\left(  bc\right)  =f
\end{multline*}

\item $e=bc=c$ and $cb=c$

We prove that $cd=c$ by showing that $\left(  a,b,cd\right)  $ is a
nonassociative triple. Indeed, $\left(  ab\right)  \left(  cd\right)
=d\left(  cd\right)  =dc=g$, while $a\left(  b\left(  cd\right)  \right)  =f$
can be derived in the following way.%
\begin{multline*}
a\left(  b\left(  cd\right)  \right)  \overset{SH}{=}a\left(  \left(
bc\right)  d\right)  =a\left(  cd\right)  \overset{SH}{=}\left(  ac\right)
d=\left(  ac\right)  \left(  ab\right) \\
\overset{SH}{=}\left(  \left(  ac\right)  a\right)  b\overset{\mathcal{B}}%
{=}\left(  ac\right)  b\overset{SH}{=}a\left(  cb\right)  =ac=a\left(
bc\right)  =f
\end{multline*}

Now we can compute that $ca=\left(  cd\right)  a=c\left(  da\right)  =c\left(
\left(  ab\right)  a\right)  =c\left(  ab\right)  =cd=c$, and the rest of the
multiplication table of $\mathbb{G}_{6}$ is not hard to fill out (we set
$h=ba$ and $i=bf$). The only entries whose calculation is not straightforward
are $ag$, $ai$ and $di$. Since $f\neq g$, at least one of these two elements
is different from $c$, hence $\left(  a,d,f\right)  $ or $\left(
a,d,g\right)  $ is an associative triple (even if $d=b$). Therefore we have
either $ag=a\left(  df\right)  =\left(  ad\right)  f=df=g$, or $ag=a\left(
dg\right)  =\left(  ad\right)  g=dg=g$ (after computing $df=dg=g$ and $ad=d$,
which is easy). Writing $ai$ either as $a\left(  bf\right)  $ or $a\left(
bg\right)  $ and $di$ as $d\left(  bf\right)  $ or $d\left(  bg\right)  $ we
get by a similar argument that $ai=g$ and $di=g$.

There are four congruences of $\mathbb{G}_{6}$ that do not collapse $f$ and
$g$, the corresponding factor groupoids are $\mathbb{G}_{7},\mathbb{G}%
_{8},\mathbb{G}_{9}$ and $\mathbb{G}_{10}$.%
\[
\mathbb{G}_{6}:\
\begin{tabular}
[c]{l|llllllll}%
$\cdot$ & $a$ & $b$ & $c$ & $d$ & $f$ & $g$ & $h$ & $i$\\\hline
$a$ & $a$ & $d$ & $f$ & $d$ & $f$ & $g$ & $d$ & $g$\\
$b$ & $h$ & $b$ & $c$ & $h$ & $i$ & $i$ & $h$ & $i$\\
$c$ & $c$ & $c$ & $c$ & $c$ & $c$ & $c$ & $c$ & $c$\\
$d$ & $d$ & $d$ & $g$ & $d$ & $g$ & $g$ & $d$ & $g$\\
$f$ & $f$ & $f$ & $f$ & $f$ & $f$ & $f$ & $f$ & $f$\\
$g$ & $g$ & $g$ & $g$ & $g$ & $g$ & $g$ & $g$ & $g$\\
$h$ & $h$ & $h$ & $i$ & $h$ & $i$ & $i$ & $h$ & $i$\\
$i$ & $i$ & $i$ & $i$ & $i$ & $i$ & $i$ & $i$ & $i$%
\end{tabular}
\]%
\[
\mathbb{G}_{7}:\
\begin{tabular}
[c]{l|lllllll}%
$\cdot$ & $a$ & $b$ & $c$ & $d$ & $f$ & $g$ & $h$\\\hline
$a$ & $a$ & $d$ & $f$ & $d$ & $f$ & $g$ & $d$\\
$b$ & $h$ & $b$ & $c$ & $h$ & $g$ & $g$ & $h$\\
$c$ & $c$ & $c$ & $c$ & $c$ & $c$ & $c$ & $c$\\
$d$ & $d$ & $d$ & $g$ & $d$ & $g$ & $g$ & $d$\\
$f$ & $f$ & $f$ & $f$ & $f$ & $f$ & $f$ & $f$\\
$g$ & $g$ & $g$ & $g$ & $g$ & $g$ & $g$ & $g$\\
$h$ & $h$ & $h$ & $g$ & $h$ & $g$ & $g$ & $h$%
\end{tabular}
\ \qquad\mathbb{G}_{8}:\quad%
\begin{tabular}
[c]{lllllll}%
$\cdot$ & \multicolumn{1}{|l}{$a$} & $b$ & $c$ & $d$ & $f$ & $g$\\\hline
$a$ & \multicolumn{1}{|l}{$a$} & $d$ & $f$ & $d$ & $f$ & $g$\\
$b$ & \multicolumn{1}{|l}{$d$} & $b$ & $c$ & $d$ & $g$ & $g$\\
$c$ & \multicolumn{1}{|l}{$c$} & $c$ & $c$ & $c$ & $c$ & $c$\\
$d$ & \multicolumn{1}{|l}{$d$} & $d$ & $g$ & $d$ & $g$ & $g$\\
$f$ & \multicolumn{1}{|l}{$f$} & $f$ & $f$ & $f$ & $f$ & $f$\\
$g$ & \multicolumn{1}{|l}{$g$} & $g$ & $g$ & $g$ & $g$ & $g$\\
&  &  &  &  &  &
\end{tabular}
\]
$\mathbb{G}_{9}:\
\begin{tabular}
[c]{l|llllll}%
$\cdot$ & $a$ & $b$ & $c$ & $d$ & $f$ & $h$\\\hline
$a$ & $a$ & $d$ & $f$ & $d$ & $f$ & $d$\\
$b$ & $h$ & $b$ & $c$ & $h$ & $h$ & $h$\\
$c$ & $c$ & $c$ & $c$ & $c$ & $c$ & $c$\\
$d$ & $d$ & $d$ & $d$ & $d$ & $d$ & $d$\\
$f$ & $f$ & $f$ & $f$ & $f$ & $f$ & $f$\\
$h$ & $h$ & $h$ & $h$ & $h$ & $h$ & $h$%
\end{tabular}
\ \qquad\mathbb{G}_{10}:\
\begin{tabular}
[c]{llllll}%
$\cdot$ & \multicolumn{1}{|l}{$a$} & $b$ & $c$ & $d$ & $f$\\\hline
$a$ & \multicolumn{1}{|l}{$a$} & $d$ & $f$ & $d$ & $f$\\
$b$ & \multicolumn{1}{|l}{$d$} & $b$ & $c$ & $d$ & $d$\\
$c$ & \multicolumn{1}{|l}{$c$} & $c$ & $c$ & $c$ & $c$\\
$d$ & \multicolumn{1}{|l}{$d$} & $d$ & $d$ & $d$ & $d$\\
$f$ & \multicolumn{1}{|l}{$f$} & $f$ & $f$ & $f$ & $f$\\
&  &  &  &  &
\end{tabular}
$
\end{caselist}

\qedhere
\end{proof}

\begin{remark}
Let us mention that there is a third possibility to measure associativity with
the help of the Hamming distance of multiplication tables. This yields the
notion of the \emph{semigroup distance} of a groupoid. Groupoids with small
semigroup distance, and connections between the semigroup distance and the
index of nonassociativity were studied in \cite{KT2 sdist}.

The different ways of measuring associativity do not seem to be closely
related. For example, the groupoid $\mathbb{G}_{3}$ is an SH-groupoid, with
the largest possible associative spectrum: $s_{\mathbb{G}_{3}}\left(
n\right)  =C_{n-1}$ for every $n$. (For the proof of the latter fact see 5.1
in \cite{CsW spec}; $\mathbb{G}_{3}$ is isomorphic to the groupoid with number
$17$ there.)

Therefore it is not surprising that the class of groupoids found in Theorem
\ref{THM spec} is disjoint from the class described in Theorem \ref{THM SH},
i.e. there is no groupoid with a minimal clone that is almost associative in
both the `spectral' and the `index' sense. Indeed, if $\mathbb{A}$ satisfies
the conditions of both theorems, then $\mathbb{A}$ (or its dual) satisfies
$x\left(  y\left(  zu\right)  \right)  \approx x\left(  \left(  yz\right)
u\right)  $ by the considerations preceding Lemma \ref{lemma general ass}, and
$\mathbb{A}$ (or its dual) contains a subgroupoid isomorphic to one of the
groupoids $\mathbb{G}_{1},\ldots,\mathbb{G}_{10}$ by Theorem
\ref{thm minimal SH in B}. However, this is impossible, because neither of
these ten groupoids and neither of their duals satisfy $x\left(  y\left(
zu\right)  \right)  \approx x\left(  \left(  yz\right)  u\right)  $ as it can
be seen from their multiplication tables (let $x=a,y=a,z=b,u=c$ for
$\mathbb{G}_{1},\ldots,\mathbb{G}_{10}$ and $x=a,y=c,z=b,u=a$ for their duals).
\end{remark}


\begin{thebibliography}{99}                                                                                               %


\bibitem {BKKR Galois}V. G. Bodnar\v{c}uk, L. A. Kalu\v{z}nin, V. N. Kotov, B.
A. Romov, \textit{Galois theory for Post algebras I-II, }Kibernetika (Kiev)
\textbf{3} (1969), 1--10; \textbf{5} (1969), 1--9. (Russian)

\bibitem {Climescu1}A. C. Climescu, \textit{\'{E}tudes sur la th\'{e}orie des
syst\`{e}mes multiplicatifs uniformes I. L'indice de non-associativit\'{e},}
Bull. \'{E}cole Polytech. Jassy \textbf{2} (1947), 347--371. (French)

\bibitem {Climescu2}A. C. Climescu, \textit{L'ind\'{e}pendance des conditions
d'associativit\'{e},} Bull. Inst. Polytech. Jassy \textbf{1} (1955), 1--9. (Romanian)

\bibitem {Cs 3all}B. Cs\'{a}k\'{a}ny, \textit{All minimal clones on the
three-element set,} Acta Cybernet. \textbf{6} (1983), no. 3, 227--238.

\bibitem {Cs cons}B. Cs\'{a}k\'{a}ny, \textit{On conservative minimal
operations, }Lectures in Universal Algebra (Szeged, 1983), Colloq. Math. Soc.
J\'{a}nos Bolyai, \textbf{43}, North-Holland, Amsterdam, 1986, 49--60.

\bibitem {CsW spec}B. Cs\'{a}k\'{a}ny, T. Waldhauser, \textit{Associative
spectra of binary operations,} Mult.-Valued Log. \textbf{5} (2000), no. 3, 175--200.

\bibitem {DrapalKepka}A. Dr\'{a}pal, T. Kepka, \textit{Sets of associative
triples,} Europ. J. Combinatorics \textbf{6} (1985), 227--231.

\bibitem {G Galois}D. Geiger, \textit{Closed systems and functions of
predicates, }Pacific J. Math. \textbf{27} (1968), 95--100.

\bibitem {Hajek1}P. H\'{a}jek, \textit{Die Sz\'{a}szschen Gruppoide,}
Mat.-Fys. \v{C}asopis Sloven. Akad. Vied \textbf{15} (1965) no. 1., 15--42. (German)

\bibitem {Hajek2}P. H\'{a}jek, \textit{Berichtigung zu meiner arbeit
\textquotedblleft Die Sz\'{a}szschen Gruppoide\textquotedblright,} Mat.-Fys.
\v{C}asopis Sloven. Akad. Vied \textbf{15} (1965) no. 4., 331. (German)

\bibitem {JQ cons}J. Je\v{z}ek, R. W. Quackenbush, \textit{Minimal clones of
conservative functions, }Internat. J. Algebra Comput. \textbf{5} (1995), no.
6, 615--630.

\bibitem {KK}K. A. Kearnes, \textit{Minimal clones with abelian
representations,} Acta Sci. Math. (Szeged) \textbf{61} (1995), no. 1-4, 59--76.

\bibitem {KSz comm}K. A. Kearnes, \'{A}. Szendrei, \textit{The classification
of commutative minimal clones,} Discuss. Math. Algebra Stochastic Methods
\textbf{19} (1999), no. 1, 147--178.

\bibitem {KT1 ns}T. Kepka, M. Trch, \textit{Groupoids and the associative law
I. (Associative triples),} Acta Univ. Carol. Math. Phys. \textbf{33} (1992),
no. 1., 69--86.

\bibitem {KT2 sdist}T. Kepka, M. Trch, \textit{Groupoids and the associative
law II. (Groupoids with small semigroup distance),} Acta Univ. Carol. Math.
Phys. \textbf{34} (1993), no. 1., 67--83.

\bibitem {KT3 SH aaa}T. Kepka, M. Trch, \textit{Groupoids and the associative
law III. ( Sz\'{a}sz-H\'{a}jek groupoids),} Acta Univ. Carol. Math. Phys.
\textbf{36} (1995), no. 1., 17--30.

\bibitem {KT4 aba}T. Kepka, M. Trch, \textit{Groupoids and the associative law
IV. (Sz\'{a}sz-H\'{a}jek groupoids of type }$\mathit{(a,b,a)}$\textit{),} Acta
Univ. Carol. Math. Phys. \textbf{35} (1994), no. 1., 31--42.

\bibitem {KT5 aab}T. Kepka, M. Trch, \textit{Groupoids and the associative law
V. (Sz\'{a}sz-H\'{a}jek groupoids of type }$\mathit{(a,a,b)}$\textit{),} Acta
Univ. Carol. Math. Phys. \textbf{36} (1995), no. 1., 31--44.

\bibitem {KT6 abc}T. Kepka, M. Trch, \textit{Groupoids and the associative law
VI. (Sz\'{a}sz-H\'{a}jek groupoids of type }$\mathit{(a,b,c)}$\textit{),} Acta
Univ. Carol. Math. Phys. \textbf{38} (1997), no. 1., 13--22.

\bibitem {LLPPP}L. L\'{e}vai, P. P. P\'{a}lfy, \textit{On binary minimal
clones,} Acta Cybernet. \textbf{12} (1996), no. 3, 279--294.

\bibitem {PPP prep}P. P. P\'{a}lfy, \textit{Minimal clones,} Preprint of the
Math. Inst. Hungarian Acad. Sci. 27/1984.

\bibitem {Plonka idred}J. P\l onka, \textit{On groups in which idempotent
reducts form a chain,} Colloq. Math. \textbf{29} (1974), 87--91.

\bibitem {Plonka k-cyc}J. P\l onka, \textit{On }$\mathit{k}$\textit{-cyclic
groupoids,} Math. Japon. \textbf{30} (1985), no. 3, 371--382.

\bibitem {Post}E. Post, \textit{The two-valued iterative systems of
mathematical logic, }Annals of Mathematics Studies, no. 5, Princeton
University Press, Princeton, 1941.

\bibitem {R 5typ}I. G. Rosenberg, \textit{Minimal clones I. The five types,}
Lectures in Universal Algebra (Szeged, 1983), Colloq. Math. Soc. J\'{a}nos
Bolyai, \textbf{43}, North-Holland, Amsterdam, 1986, 405--427.

\bibitem {enc}N. J. A. Sloane, \textit{The On-Line Encyclopedia of Integer
Sequences,}\newline http://www.research.att.com/$\sim$njas/sequences, 2005.

\bibitem {Szasz}G. Sz\'{a}sz, \textit{Die Unabh\"{a}ngigkeit der
Assoziativit\"{a}tsbedingungen,} Acta Sci. Math. (Szeged) \textbf{15} (1953),
20--28. (German)

\bibitem {Szcz}B. Szczepara, \textit{Minimal clones generated by groupoids,
}Ph.D. Thesis, Universit\'{e} de Montr\'{e}al, 1995.

\bibitem {SzA clUA}\'{A}. Szendrei, \textit{Clones in Universal Algebra,}
S\'{e}minaire de Math\'{e}matiques Sup\'{e}rieures, \textbf{99}, Presses de
L'Universit\'{e} de Montr\'{e}al, 1986.

\bibitem {SzM minfcs}M. B. Szendrei \textit{On closed sets of term functions
on bands,} Semigroups (Proc. Conf., Math. Res. Inst., Oberwolfach, 1978), pp.
156--181, Lecture Notes in Math., 855, Springer, Berlin, 1981.

\bibitem {W 4maj}T. Waldhauser, \textit{Minimal clones generated by majority
operations,} Algebra Universalis \textbf{44} (2000), no. 1-2, 15--26.

\bibitem {W wa}T. Waldhauser, \textit{Minimal clones with weakly abelian
representations, }Acta Sci. Math. (Szeged) \textbf{69} (2003), no. 3-4, 505--521.
\end{thebibliography}
\end{document}